\documentclass[a4paper]{amsart}

\usepackage{amsfonts,amscd,latexsym,amssymb}
\usepackage{eucal,mathrsfs}
\usepackage[all]{xy}
\usepackage{pxfonts}

\newcommand{\B}[1]{{\mathbb #1}}
\newcommand{\C}[1]{{\mathcal #1}}

\font\tencyr=wncyr10
\font\sevencyr=wncyr7
\font\fivecyr=wncyr5
\newfam\cyrfam
\textfont\cyrfam=\tencyr
\scriptfont\cyrfam=\sevencyr
\scriptscriptfont\cyrfam=\fivecyr
\def\bukfa{\fam\cyrfam\tencyr}

\font\tencyr=wncyr10
\font\sevencyr=wncyr7
\font\fivecyr=wncyr5
\newfam\cyrfam
\textfont\cyrfam=\tencyr
\scriptfont\cyrfam=\sevencyr
\scriptscriptfont\cyrfam=\fivecyr

\newtheorem{theorem}[subsection]{Theorem}

\newtheorem{corollary}[subsection]{Corollary}
\newtheorem{lemma}[subsection]{Lemma}
\newtheorem{proposition}[subsection]{Proposition}
\theoremstyle{definition}
\newtheorem{definition}[subsection]{Definition}
\newtheorem{example}[subsection]{Example}
\theoremstyle{remark}
\newtheorem{remark}[subsection]{Remark}

\numberwithin{figure}{section}
\numberwithin{table}{section}

\newcommand{\al}{{\alpha}}

\newcommand{\be}{{\beta}}

\newcommand{\Om}{{\Omega}}
\newcommand{\om}{{\omega}}

\newcommand{\ga}{{\gamma}}
\newcommand{\Ga}{{\Gamma}}

\newcommand{\la}{{\lambda}}

\newcommand{\si}{{\sigma}}
\newcommand{\Si}{{\Sigma}}
\newcommand{\vfi}{{\varphi}}
\newcommand{\cp}{{\B C\B P}}
\newcommand{\pci}{{\B C\B P^{\infty}}}

\newcommand{\Mo}{(M,\omega )}
\newcommand{\Mwo}{(M,\widehat \omega )}
\newcommand{\Wo}{(W,\omega _W)}

\newcommand\BS{\operatorname{BSymp}}
\newcommand\BH{\operatorname{BHam}}

\newcommand\Flux{\operatorname{Flux}}
\newcommand\Gr{\operatorname{Gr}}
\newcommand\Ham{\operatorname{Ham}}

\newcommand\Hom{\operatorname{Hom}}

\newcommand\Map{\operatorname{Map}}

\newcommand\Symp{\operatorname{Symp}}

\newcommand{\bs}{{\bigskip}}
\newcommand{\NI}{{\noindent}}

\begin{document}

\title{Symplectic configurations}
\author{\'Swiat Gal}
\address{Gal:\newline
Mathematical Institute \\
University of Wroc\l aw \\
pl. Grunwaldzki 2/4 \\
50-384 Wroc\l aw \\
Poland\newline 
Institute of Mathematics \\
University of Warsaw \\
ul. Banacha 2 \\
02-097 Warszawa \\
Poland}
\email{sgal@math.uni.wroc.pl}
\urladdr{http://www.math.uni.wroc.pl/\~{}sgal}
\author{Jarek K\c edra}
\address{K\c edra:\newline
Institute of Mathematics \\
University of Szczecin \\
Wielkopolska 15 \\
70-451 Szczecin \\
Poland \newline
Department of Mathematical Sciences \\
University of Aberdeen \\
Meston Building \\
King's College \\
Aberdeen AB24 3UE \\
Scotland \\
UK}
\email{kedra@maths.abdn.ac.uk}
\urladdr{http://www.maths.abdn.ac.uk/\~{}kedra}
\date{\today }
\thanks{The first author is partially supported by the State Committee
for Scientific Research (grant no. 2 P03A 017 25) and by the Foundation
for Polish Science.
\newline\indent
The second author is partially supported by the State Committee for Scientific
Research (grant no. 1 PO3A 023 27).} 


\begin{abstract}
\textfont\cyrfam=\sevencyr
\scriptfont\cyrfam=\fivecyr
We define a class of symplectic fibrations called symplectic
configurations. They are a natural generalization of Hamiltonian
fibrations in the sense that they
admit a closed symplectic connection two--form.
Their geometric and topological properties are 
investigated. We are mainly concentrated on integral symplectic
manifolds.

We construct the classifying space $\bukfa B$ of
symplectic integral configurations. The properties of
the classifying map ${\bukfa B}\to \BS\Mo$ are examined.
The universal symplectic bundle over $\bukfa B$ has
a natural connection whose holonomy group is 
the enlarged Hamiltonian group recently defined by
McDuff.
The space $\bukfa B$ is identified with the classifying
space of a certain extension of 
the symplectomorphism group. 

\end{abstract}


\maketitle
{
\setcounter{tocdepth}1
\tableofcontents
}

\section{Motivation}\label{S:motivic}
The paper introduces a new object, called a {\em symplectic configuration},
which is a symplectic fibration equipped with a map of the total space
into some symplectic manifold so that the restrictions of the map to the
fibers are symplectic embeddings. Symplectic configurations provide
a generalization of Hamiltonian fibrations in the sense that they
admit a closed symplectic connection two--form
\cite[Theorem 6.21]{MR2000g:53098}.

We are partially motivated by the result of
Narasimhan and Ramanan \cite{MR0133772}. They 
showed that the natural connection $\la$
on the universal principal bundle
$$U(n)\to E \to \Gr(n,\infty)$$
has the following universal property.
For any principal $U(n)$-bundle $P\to B$ equipped 
with a connection $\lambda_P$, there
exists a map $f\colon B\to \Gr(n,\infty)$ 
such that $(P,\lambda_P)=f^*(E,\lambda)$.
Notice that the total space $E$ is the space of 
linear embeddings from $\B C^n$ into $\B C^\infty$
and $\Gr(n,\infty)$ is the Grassmannian 
of $n$-dimensional subspaces in~$\B C^{\infty}$.

This paper is devoted to proving a similar result about the group of
symplectomorphisms $\Symp\Mo$ in place of $U(n)$.
Mostly, we deal with the case when $\om$ has integral periods.
Assume for the moment that $H^2(M,\B Z)$ is torsion free (the torsion issues
are discussed in detail in Section \ref{S:intro}).

The space $\Symp(M,\pci)$ of symplectic embeddings of $M$ into $\pci$  
is not contractible (it has non-trivial $\pi_2$, see Corollary \ref{C:vanish}).
Thus it is not the total space of the universal $\Symp\Mo$-bundle. 
However, we find a deep connection with the 
result of Narasimhan and Ramanan mentioned above. 
That is, we construct
a natural connection on the principal bundle
$$
\Symp\Mo\to\Symp(M,\pci)\to\Symp(M,\pci)/\Symp\Mo=:{\bukfa B}_M.
$$

The connections on  $\Symp\Mo$-bundles are in one--to-one correspondence
with vertically closed extensions of $\omega$ to the total space of the
associated bundle with fiber $M$ \cite[Lemma 6.18]{MR2000g:53098}. 
Among them, the connections represented
by closed two--forms are of special interest since they
generalize the notion of a coupling form introduced by
Guillemin, Lerman and Sternberg
\cite{MR98d:58074}.
A coupling form is a certain closed connection two--form
satisfying a normalization condition.
Furthermore, the fiber integrals of powers of the cohomology class
of the coupling form give a sequence of symplectic characteristic
classes \cite{JK,MR2115670}.

The connection two--form $\Om$ that we construct in the 
configuration 
$\Mo\to M_{\bukfa B}\to {\bukfa B}_M$ 
is universal in the following sense.
For any symplectic bundle $M\to E\to B$ endowed
with a closed connection two--form $\Omega_E$ there exists a map
$f\colon B\to{\bukfa B}_M$ such that $(E,\Om_E)=f^*(M_{\bukfa B},\Om)$.

We show that ${\bukfa B}_M$ is a classifying space of 
an extension of $\Symp\Mo$ by the gauge group $\Map(M,U(1))$.
We determine the holonomy subgroup of
this connection and prove that it is
equal to the subgroup $\Ham^{s\B Z}\Mo$ of 
$\Symp\Mo$ recently discovered by McDuff \cite{dusa}.
Her paper  deals with , among other things, a 
characterization of symplectic bundles 
with closed connection form.
It was important to the development of the present paper 
that Dusa McDuff shared with us her ideas on a generalization of 
Hamiltonian fibrations.

\subsection*{Organization of the paper}

In Section \ref{S:intro} we give the necessary definitions and state
the main results. We give the direct parts of the proofs and postpone
the more technical ones to later sections.

In Section \ref{S:top} we prove the universality of the 
configuration $\Mo\to M_{\bukfa B}\to \bukfa B$
and prove the homotopy properties of the space $\bukfa B$.
We also prove that $\bukfa B$ is the classifying space
of an extension $\C G$ of $\Symp\Mwo$.

Section \ref{S:gsc} is devoted to the geometry of symplectic
configurations. More precisely, we define a principal connection
on a fibration $\Symp(M,W)\to {\bukfa B}_W$ and investigate its
properties.

In Section \ref{S:group} we investigate the group cohomology
relations between groups $\C D$ and $\C G$.
We also compare these groups with the McDuff 
subgroup $\Ham^{s\B Z}\Mo$.
The section starts with some preparation on crossed homomorphisms
and the necessary algebraic topology.

In Section \ref{S:cohomology} we investigate characteristic classes
of integral symplectic configurations and relate them to other known
characteristic classes.

\subsection*{Acknowledgement}
The idea of configurations  grew out of discussions
with Tadeusz Januszkiewicz. 

We warmly thank Dusa McDuff for showing us an early version
of her paper \cite{dusa}, discussions, comments and
picking up mistakes and inaccuracies in earlier drafts 
of this paper.
In particular, she drew our attention to the torsion
issues.

We thank Agnieszka Bojanowska, Stefan Jackowski and
Michael Weiss for homotopy discussions and hospitality.

\section{Preliminaries and statements of results}\label{S:intro}

\subsection{Symplectic configurations}
Let $\Mo$ be a closed symplectic manifold.
A symplectic manifold is called {\em integral} if the 
symplectic form has integral periods.

\begin{definition}\label{D:configuration}
Let $\Mo$ and $\Wo$ be symplectic manifolds. We say that a symplectic
fiber bundle (a fiber bundle with a structure group $\Symp\Mo$)
$M\to E\to B$ is a 
{\em $W$-symplectic configuration} 
if there exists a map $E\to W$ whose restriction
to any fiber of $E$ is a symplectic embedding.
The fibration $E$ is called an (integral)  symplectic configuration
if it is a $W$-configuration for some (integral) $\Wo$.
\end{definition}

Here is an alternative approach.
Consider the space $\Symp(M,W)$ of symplectic embeddings
$f\colon \Mo\to \Wo$.
The group of all symplectomorphisms of the source acts freely on that
space and the quotient is denoted by ${\bukfa B}_W$.
We call it the {\em space of symplectic configurations
of $\Mo$ in $\Wo$}. We get a principal fibration
$$\Symp\Mo \to \Symp(M,W)\to {\bukfa B}_W$$
and the associated symplectic one
$$\Mo\to M_W\to {\bukfa B}_W.$$

\begin{proposition}\label{P:def}
A symplectic fibration $\Mo\to E\to B$ is a $W$-symplectic configuration 
if and only if it is a pull-back of the bundle $M_W$.
\end{proposition}
\begin{proof}
Let  $E\to W$ be a map such that it is a symplectic
embedding on every fiber.
Clearly, it defines a map $B\to {\bukfa B}_W$
such that $E$ is a pull-back of $M_W$.

To prove the converse
choose a point $p\in M$ and
consider the following composition of two maps:
$$
\CD
\Symp(M,W) @>\C P>> M_W @>ev>> W,
\endCD
$$
where
$$
\C P(f):= [f,p]\quad {\text {and}}\quad
ev[f,q]=f(q).
$$
The map $ev$ is well defined.
Indeed $M_W$ is the quotient of $\Symp(M,W)\times M$
by the relation $[f\circ\psi,q]\sim[f,\psi(q)]$
and $ev([f\circ \psi,q])=f(\psi(q))=ev([f,\psi(q)])$.
The evaluation map is the required map for the fibration
$M_W$. Since $E$ is a pull-back of $M_W$, the proof is finished.
\end{proof}

\subsection{Closed connection two--forms}
Without loss of generality we will assume that the base of a fibration
is a manifold, so we may speak of the differential forms on~$E$.

Recall that if $\Mo\to E\to B$ is a symplectic fibration
then a closed two--form
$\Om\in \Om^2(E)$ is called a {\em closed symplectic connection
two--form} if it restricts to the symplectic form on any fiber
(see Lemma 6.18 in \cite{MR2000g:53098} for a general definition).
According to Thurston (Theorem 6.3 in 
\cite{MR2000g:53098}), the existence of a closed connection two--form
is a purely cohomological condition,
i.e. the form $\Om$ exists if and only if the class $[\om]$
belongs to the image of the map 
$H^2(E;\B R)\to H^2(M;\B R)$.

\begin{lemma}\label{L:two--form}
A symplectic bundle $\Mo\to E\to B$ 
is a symplectic configuration
if and only if
it admits a closed connection two--form.
\end{lemma}

\begin{proof} If $f\colon E\to W$ defines a configuration, then 
$f^*(\om_W)$ is a closed connection two--form on $E$.

To prove the converse construct a map
$f\colon(E,\Om)\to (W,\om_W):=\varprod_i(\cp^{n_i}, a_i \om_{n_i})$
such that $f^*[\om_W]=[\Om]$ for appropriate $a_i\in \B R$, 
where $\om_n$ is the standard symplectic form on $\cp^n$.
For each $\eta_i\in H^2(E,\B Z)$ there exists a map
$f_i\colon E \to \cp^{\infty}$ such that $f_i^*[\om_{\infty}] = \eta_i$
(here $[\omega_{\infty}]$ is the generator
of $H^2(\cp^{\infty};\B Z)$).
Hence $[\Om]=\sum a_i\eta_i =
(\varprod f_i)^*(a_1[\om_{\infty}]\times \dots \times a_n[\om_{\infty}])$
for the product map $\varprod f_i\colon E \to \varprod_i \cp^{\infty}$.
The latter factors through a product
of finite dimensional projective spaces, as required.
For sufficiently large
$\sum n_i$ we can perturb $f$ into an
embedding preserving the two--form, according to the
h-principle (see Section 3.4.2 in \cite{MR90a:58201}).
\end{proof}

\subsection{The universal integral configuration}
Due to Gromov, every  integral symplectic manifold
embeds into complex projective space $(\cp^N,\om_N)$ equipped
with the standard symplectic form for $N$ big enough 
(Exercise 3.4.2.(1) in \cite{MR90a:58201}). 
Spaces of symplectic embeddings satisfy the obvious 
functoriality properties and we define
$\Symp(M,\pci):=\bigcup_{N\to \infty} \Symp(M,\cp^N)$.
Its quotient by $\Symp\Mo$ is denoted by ${\bukfa B}_{\infty}$ and called
{\em the universal symplectic configuration space}.
As before we have the universal principal fibration
$$
\Symp\Mo\to \Symp(M,\pci) \to {\bukfa B}_{\infty}
$$
and the associated one
$$
\Mo \to M_{{\bukfa B}_{\infty}} \to {\bukfa B}_{\infty}.
$$
Let  $\Mo\to E\to B$ be a symplectic fibration
whose classifying map $\xi\colon B\to \BS\Mo$ has a lift
$\widehat \xi\colon B\to \bukfa B_{\infty}$.
This lift defines a map
$E\to M_{{\bukfa B}_{\infty}}\stackrel{ev}\to \pci$
which factors through $\cp^N$, since $E$ is finite dimensional.
Thus $E$ is an integral configuration
(see Definition \ref{D:configuration}),
and the pull-back of the standard symplectic form
$\om_N$ on $\cp^N$ is a closed integral connection form on $E$.


\subsection{Torsion}\label{SS:torsion} 
The space $\Symp(M,\pci)$ has the same number of connected components
as the set of maps from M to $\pci$ with the property that the generator of 
$H^2(\pci,\B R)$ pulls back to $[\omega]$.
The homotopy classes of maps from $M$ to $\pci$
are in one--to-one correspondence with $H^2(M,\B Z)$.
The kernel of the map $H^2(M,\B Z)\to
H^2(M,\B R)$ consists of the torsion of $H^2(M,\B Z)$.
Thus if this torsion in nontrivial then
the space $\Symp(M,\pci)$ is disconnected.
Since $M$ is closed, $\Symp(M,\pci)$ has a finite
number of connected components.
Moreover, if $\Symp(M,\omega)$ acts non-transitively
on the preimage of $\omega$ in $H^2(M,\B Z)$
then $\bukfa B$ is disconnected.

Let us observe that any (finite) Abelian group $T$
can be the torsion of $H^2(M;\B Z)$. First notice that
the torsion of $H^2(M;\B Z)$ is isomorphic to the
torsion of $H_1(M;\B Z)$ (see Corollary 5.5.4 in \cite{MR35:1007}).
Take a finitely presented group $G$ such that its abelianization
has torsion equal to $T$ (e.g. $G=T$). Due to Gompf \cite{MR96j:57025}, 
there exists a closed symplectic 4-manifold $\Mo$ with 
fundamental group $G$. Hence the torsion of $H^2(M;\B Z)$
is isomorphic to $T$.

Let $\Mo$ be an integral symplectic manifold and let
$\widehat \omega$ denote a pair $(\om,\tau)$ of the
symplectic form and a lift of its class to
$H^2(M,\B Z)$. Define 
$$\Symp\Mwo:=
\{\psi\in \Symp\Mo\colon  \psi^*(\tau)=\tau\}$$
to be a finite index open-closed
subgroup of $\Symp\Mo$
consisting of the symplectomorphisms which
preserve the class $\tau$. 
Notice that if $H^2(M,\B Z)$ is 
torsion-free then $\Symp\Mwo =\Symp\Mo$.

\begin{lemma}
If $\Mo\to E\to B$ is an integral symplectic configuration,
then there exists a lift 
$\tau\in H^2(M,\B Z)$ of $[\om]$ such that the 
structure group of $E$ reduces to $\Symp\Mwo$.
\end{lemma}
\begin{proof}
If $f\colon E\to\pci$ defines
a configuration structure on $E$ then the
pull-back of the generator of $H^2(\pci,\B Z)$
has the above property.
\end{proof}
Let $\Symp_\tau(M,\pci)$ denote the connected component
of $\Symp(M,\pci)$ containing symplectic embeddings corresponding
to $\tau\in H^2(M,\B Z)$. The quotient 
${\bukfa B}={\bukfa B}_\tau=\Symp_\tau(M,\pci) /\Symp\Mwo$ 
is the component of ${\bukfa B}_{\infty}$ corresponding to $\tau$.
In the sequel, we restrict ourselves to the investigation of
${\bukfa B}$. We will also call it a universal configuration
space. Again we have the universal fibration:
$$
\Symp\Mwo \to \Symp_\tau(M,\pci) \to \bukfa B
$$
and the associated one
$$
\Mo \to M_{\bukfa B} \to \bukfa B.
$$
A symplectic fibration whose classifying map
lifts to ${\bukfa B}_{\tau}$ will be called an integral symplectic
configuration compatible with $\tau$. Notice that
a symplectic integral configuration might be compatible
with several lifts.

Since $\Symp\Mwo$ is open-closed in $\Symp\Mo$, their
connected components of the identity are equal. Therefore
the map of classifying spaces induced by the inclusion is
a  covering. In particular, the higher homotopy groups of these
classifying spaces are equal, as well as  rational 
(co-) homologies. 

\subsection{Homotopy properties of the universal configuration space}
\label{SS:BS}
Integral symplectic configurations constitute
a big class of symplectic fibrations. To see this
consider the classifying map
$$
\C F\colon {\bukfa B}\to \BS\Mo
$$
and the maps
$\C F_k\colon \pi_k({\bukfa B})\to \pi_k(\BS\Mo) $ 
induced on homotopy groups. 
The following theorem is
proved in Section \ref{SS:TBS}.

\begin{theorem}\label{T:BS} 
Let $\Mo $ be a compact integral symplectic manifold. Then
\begin{enumerate}
\item
$\pi_k({\bukfa B}) \stackrel{\C F_k}\cong \pi_k(\BS\Mo)$ for $k\neq 1,2.$
\item
The following sequences are exact.
\begin{eqnarray*}
0\to H^0(M;\B Z)\to \pi_2({\bukfa B}) 
&\stackrel{\C F_2}\to& \pi_2(\BS \Mo) \to \Gamma_{\om}\to 0 \\
0\to \Gamma_{\om}\to H^1(M;\B Z)\to \pi_1({\bukfa B})
&\stackrel{\C F_1}\to& \pi_1(\BS\Mwo)\to 0,
\end{eqnarray*}
where $\Gamma_\omega$ is the flux group (see \cite[p.~321]{MR2000g:53098}).
\end{enumerate}
\end{theorem}

\subsection{An extension of $\Symp\Mwo$}
We identify the universal symplectic configuration space $\bukfa B$
as the classifying space $B\C G_{\tau}$ of a certain extension $\C G_{\tau}$
of $\Symp\Mwo$. More precisely this extension is of the form
$$
0\to \Map(M,U(1))\to \C G_{\tau}\to \Symp\Mwo \to 1,
$$
where $\Map(M,U(1))$ is the (Abelian) gauge group of continuous maps
from $M$ to $U(1)$ (see Section \ref{SS:gg} for the precise definition).
We shall usually omit the subscript $\tau$ if it does not 
lead to  confusion.
The proof of the next theorem is in Section \ref{SS:gg}.

\begin{theorem}\label{T:gg}
The classifying space of $\C G_{\tau}$ is homotopy equivalent to the universal
symplectic configuration space:
$$
B\C G_{\tau} = \bukfa B_{\tau}.
$$
Moreover, the classifying map 
$\C F\colon {\bukfa B}_{\tau}=B{\C G}_{\tau}\to \BS\Mwo$
is induced by the projection ${\C G}_{\tau}\to\Symp\Mwo$.
\end{theorem}

\subsection{The holonomy group}
If $\Om$ is the symplectic form on $\pci$ and $ev\colon M_{\bukfa B}\to\pci$
is the evaluation map, then $ev^*\Om$ is the closed connection form on
the universal symplectic configuration $\Mo\to M_{\bukfa B}\to {\bukfa B}$.
It defines a symplectic connection whose holonomy group
is denoted by $\C D=\C D_{\tau}\subset \Symp\Mwo$. 

Note that $\C D$ depends on the lift $\tau$ and is defined only
up to conjugacy. In general we will omit the subscript $\tau$ except
in necessary cases.

\begin{theorem}\label{T:gd}
\begin{enumerate}
\item
$\C D$ intersects every connected component
of $\Symp\Mwo$;
\item
The identity component $\C D_0$ of $\C D$ is
equal to $\Ham \Mo$;
\item
$\C D \cap \Symp_0\Mo = \Flux^{-1}(H^1(M;\B Z)/\Ga_{\om})$, that is,
the intersection consists of symplectomorphisms with an integral
flux. In particular, $\C D$ is a closed subgroup of
$\Symp\Mwo$.
\end{enumerate}
\end{theorem}
\begin{proof}
(1) Since $\pi_1({\bukfa B})\to \pi_1(B\Symp\Mwo)$ is
a surjection (Theorem \ref{T:BS}),
$\C D$ intersects every connected component of $\Symp\Mwo$.

(2) Let $\Mo \to M_W\to {\bukfa B}_W$ be a symplectic configuration
fibration. Since $\Om_W$ is a closed connection two--form,
the holonomy about any contractible loop is Hamiltonian
due to Theorem 6.21 in \cite{MR2000g:53098}. This proves that
$\C D_0\subset \Ham\Mo$.
The converse inclusion follows from Corollary \ref{C:curv}.

(3) This part will be proved in Proposition \ref{P:overline}(2).
\end{proof}

Finally we find a connection between the holonomy subgroup $\C D$  and the
extension~$\C G$. 
\begin{theorem}\label{T:DGMap}
Let $\Map_0(M,U(1))$ denote  the identity component of $\Map(M,U(1))$.
There exists a commutative triangle of continuous homomorphisms 
$$
\xymatrix{
&\C G/\Map_0(M,U(1))\ar[d]\\
\C D\ar[r]\ar[ru]&\Symp\Mwo
}
$$
in which the diagonal arrow is a homotopy equivalence. 
In particular, the structure group
of $E$ lifts to $\C G/\Map_0(M,U(1))$ if and only if
the structure group of $E$ reduces to $\C D$.
\end{theorem}

\begin{proof}
This follows straightforward from Theorem \ref{T:dusa} and
Corollary \ref{C:section}.
\end{proof}

The next two theorems summarize all the main results of the paper.

\begin{theorem}\label{T:main}
Let $\Mo $ be a closed integral symplectic manifold and
let $\tau$ be an integral lift of $[\om]$.
Let $\Mo\to E \to B$ be a symplectic fibration.
The following conditions are equivalent:
\begin{enumerate}
\item
$E$ admits a closed integral connection two--form $\Om$ compatible with
$\tau$, i.e. there exists $\C T\in H^2(E,\B Z)$ such that it lifts
$\Om$ and extends $\tau$;
\item
$E$ is an integral symplectic configuration compatible with $\tau$;
\item
the structure group of $E$ lifts to $\C G_{\tau}$.
\end{enumerate}
\end{theorem}
\begin{proof}

\NI
{\bf (1) $\Rightarrow $ (2)} is proved in Section \ref{SS:1=>2},

\NI
{\bf (2) $\Rightarrow $ (1)}
follows from  Lemma \ref{L:two--form},

\NI
{\bf (1) $\Rightarrow $ (3)} is proved in Section \ref{SS:1=>3},

\NI
{\bf (2) $\Leftrightarrow $ (3)} follows form Theorem \ref{T:gg}.
\end{proof}

Clearly, if $E$ is an integral symplectic configuration with respect
to
the lift $\tau$ then its structure
group reduces to the holonomy group $\C D_{\tau}$. In other words,
if the structure group lifts to $\C G_{\tau}$ then it lifts
to $\C G_{\tau}/\Map_0(M,U(1))$.
The following example due to McDuff shows that the converse is not
true.

\begin{example}
Let $S^2 \to BSO(2) \stackrel{\pi}\to BSO(3)$ be the universal
fibration associated with the action of $SO(3)$ on $S^2$.
Let $\tau$ be the generator of $H^2(S^2,\B Z)$.
We claim that $\tau$ does not admit an integral extension 
although the fibration is Hamiltonian
(i.e. its structure group $SO(3)$ is a subgroup of $\Ham(S^2,area)$).
Indeed, $H^3(BSO(2);\B Z)=0$ and therefore
$H^3(BSO(3),\B Z)=\B Z/2$ is in the kernel of $\pi^*$.
Hence the generator of the latter group equals $d_2(\tau)$,
where $d_2\colon H^2(S^2) \to H^3(BSO(3))$ is the differential in the
spectral sequence. Thus $\tau$ is not in the image of 
$\pi^*\colon H^2(BSO(2))\to H^2(S^2,\B Z)$
i.e. $\tau$ does not admit an integral extension.
\end{example}

The proof of the following theorem will be given in Section
\ref{SS:4=>3}.

\begin{theorem}\label{T:ugly}
Suppose that the structure group of a symplectic fibration 
$M\to E\to B$ reduces to $\C D_{\tau}\subset \Symp\Mwo$. 
If $H_2(B,\B Z)$ is torsion-free
then $E$ is an integral configuration.
\end{theorem}

\begin{remark}
We do not prove that  $E$ is a configuration compatible with $\tau$.
We show that $d_2(\tau)=0$ and $d_3(\tau)\in E_3^{3,0}$ 
is torsion in the spectral sequence associated with $M\to E\to B$.
Thus, {\it a priori}, $\tau$  might not extend.
Notice that, although $E^{3,0}_2=H^3(B,\B Z)$ is
(by the assumption) torsion free,
$E^{3,0}_3 = E^{3,0}_2 / d_2(E^{1,1}_2)$ might have nontrivial torsion.
However, we do not know any example in which $E$ is a configuration
not compatible with $\tau$.
\end{remark}

\section{The topology of symplectic configurations}\label{S:top}

\subsection{Gromov's h-principle}\label{SS:h-principle}
Let $\Map(M,\cp^{\infty})$ (respectively $\Map^\infty(M,\pci)$ )
denotes the space of all continuous (respectively smooth) maps
from $M$ to $\cp^{\infty}$ equipped with compact-open 
(respectively $C^{\infty}$) topology. 
The space of smooth maps is defined as
$\bigcup_n C^{\infty}(M,\cp^n)$.
It is well known that the inclusion
$\Map^\infty(M,\cp^{\infty})\to \Map(M,\cp^{\infty})$
induces a (weak) homotopy equivalence. Due to a
big codimension the first space is homotopy equivalent
to the space of embeddings of $M$ into $\pci$.

\begin{theorem}\label{T:isomorphism}
Let $\tau\in H^2(M,\B Z)$ be a lift of $[\om]\in H^2(M,\B R)$.
Let $\Symp_{\tau}(M,\pci)$ and $\Map_{\tau}(M,\pci)$ denote
the connected components corresponding to the class $\tau$. 
The inclusion $i \colon \Symp_{\tau}(M,\pci) \to \Map_{\tau}(M,\pci)$ induces
an isomorphisms on homotopy groups.
\end{theorem}

\begin{remark}
The above theorem easily follows from
the parametric version of {\em h-principle for symplectic embeddings}.
Unfortunately, it seems that there is no proof of it in the
literature yet. That is why we prove the above theorem appealing only
to the well known Gromov's {\em h-principle} for symplectic immersions 
and its parametric version (Theorem (A) and Exercise (2) in 
Section 3.4.2 of \cite{MR90a:58201}).
\end{remark}

 


\begin{proof}[Proof of Theorem \ref{T:isomorphism}]

\NI
{\bf Injectivity:}
Suppose that $f\in \ker i_*$ that is there exist a commutative
diagram
$$
\xymatrix
{
S^k \times M \ar[r]^f \ar[d] & \cp^m\\
D^{k+1}\times M \ar[ru]^{F}
}
$$
where $f$ is a symplectic embedding on each fiber $\{s\}\times M$
and $F$ is a smooth embedding on fibers and equals $f$ over
the boundary. The argument consists of several steps.

We first deform $F$ to a fiberwise symplectic immersion $F'$
so that it equals $f$ over the boundary. This can be done 
according to the parametric h-principle for symplectic immersions
(Theorem 16.4.3 in \cite{MR2003g:53164}). We want to improve it
so that it will be fiberwise symplectic embedding.

We will need an isotropic embedding $j\colon M\to (\B C^n,\om_n)$.
For example one can compose an embedding of $M$ into $\B R^N$
with the standard inclusion of Lagrangian $\B R^n$ into $\B C^n$

Secondly, we need a symplectic embedding $\vfi$ of $\cp^m\times D^{2n}$ with
the product form into $\cp^N$ such that $\vfi$ is linear on $\cp^m\times {0}$.
This can be achieved by realizing $\vfi$ of $\cp^m\times D^{2n}$ in
$\cp^m\times\cp^n$ and embedding the latter in $\cp^{mn+m+n}$ via the
Segre embedding.

Take
$$
D^{k+1} \times M \to \cp^m\times D^{2n}\to \cp^N$$
defined by
$$
(d,m)\mapsto \vfi(F'(d,m),\al(d)j(m)),
$$
where $\al\colon D^{k+1}\to \B R$ is a sufficiently small
scaling function such that
it equals 0 exactly on the boundary sphere. Clearly, it is an embedding
because $j$ is. It is fiberwise symplectic because $F'$ is fiberwise
symplectic and $j$ is isotropic.

\NI
{\bf Surjectivity:}
We shall show that every fiberwise smooth embedding
$f\colon S^k\times M \to \cp^m$ (fiberwise in $\Map_{\tau}(M,\cp^m)$)
is homotopic to a fiberwise
symplectic embedding. First take the composition 
$F\colon T^*S^k\times M \to \cp^m$ of
$f$ and the projection $T^*S^k\times M\to S^k\times M$.
According to the standard h-principle there exist an arbitrarily
$C^0$-small deformation of $F$ that is a symplectic embedding
(with respect to the product symplectic form).
The restriction of it to the image of the ``zero section''
$S^k\times M\subset T^*S^k\times M$ gives a
symplectic embedding. 

The above homotopy might not
preserve the base-point but this is not a problem 
due to the injectivity of $i_*$ on $\pi_0$.
Let $s_0\in S^k$ and $m\in \Symp_{\tau}(M,\pci)$ be base-points.
Let $F\colon [0,1]\times S^k \to \Map_{\tau}(M,\pci)$ 
be the free homotopy ($F_0(s_0)=f(s_0)=m$) such that 
$F_1(S^k) \subset \Symp_{\tau}(M,\pci)$.
Let $q\colon [0,1]\to \Map_{\tau}(M,\pci)$ be defined
by $q(t):= F_t(s_0)$.

If $l\colon [0,1]\to X$ is a path in a topological space 
then $l*\colon \pi_k(X,l(0))\to \pi_k(X,l(1))$ denotes
the induced isomorphism. In our situation,
$F$ is a based homotopy from $F_0$ to $q*F_1$.
Due to the injectivity of $\pi_0(i)$, there exists
a path $p\colon [0,1]\to \Symp_{\tau}(M,\pci)$ from
$F_0(s_0)$ to $F_1(s_0)$. 
Further, we have that
$$
p*F_1 \sim (p*q^{-1})*q*F_1 \sim (p*q^{-1})*F_0\sim F_0.
$$
The last equivalence follows from the fact that
$\Map_{\tau}(M,\pci)$ is an H-space, so the action
of the fundamental group on all homotopy groups is trivial.
This finishes the proof.

\end{proof}

\begin{remark}
Note that the result we mention at the beginning of the Section \ref{S:motivic}
implies injectivity on $\pi_0$. Indeed, an integral two--form is a curvature
form of some $U(1)$-connection (on the pull-back of the universal
$U(1)$-bundle $U(1)\to S^\infty\to\pci$ by the map $f\colon M\to\pci$)
thus, following Narasimhan and Ramanan \cite{MR0133772} one can approximate $f$
(by a map classifying the same bundle, thus homotopic) to the map
preserving not only curvature $\om$ but also the connection.
\end{remark}

\subsection{Proof of Theorem \ref{T:BS}}\label{SS:TBS}

\begin{lemma}\label{L:spanier}
$\pi_k \Map(M,\pci) = H^{2-k}(M,\B Z)$, where we set
$H^m(M,\B Z)=0$ for
$m<0$.
\end{lemma}

\begin{proof}
If follows from the basic algebraic topology 
(e.g. \cite{MR35:1007}). 
Let $s\in S^k$ and $f\in \Map(M,X)$ be a base point,
the constant map $m\mapsto x\in X$.
Moreover let $M_+$ denote $M$ with an artificially added base point.
Then we calculate
\begin{eqnarray*}
\pi_k(\Map(M,X))&=&[(S^k,s),(\Map(M,X),f)]\\
&=&[M_+\wedge S^k,(X,x)]\\
&=&[M_+,\Omega^k (X,x)]
\end{eqnarray*}
Since $\cp^\infty=K(\B Z,2)$, we have that
$\pi_k(\Map(M,\pci))=[M,K(\B Z,2-k)]=H^{2-k}(M,\B Z)$.

It is easy to see a composition of
this isomorphism with the projection
$H^{2-k}(M,\B Z) \to \Hom(H_{2-k}(M,\B Z),\B Z)$ 
(which is an isomorphism for $k\in \{1,2\}$) as follows.
Let $[\xi] \in \pi_k(\Map(M,\pci))$ be represented by a map
$\xi\colon S^k\times M\to \pci $. The value of the resulting cohomology
class on $(2-k)-$cycle $a$ is equal to 
$\left <[S^k\times a],\xi^*\om_{\infty}\right >$, where 
$\om_{\infty}\in H^2(\pci ;\B Z)$ is the generator.

\end{proof}




\begin{corollary}\label{C:vanish}
$\pi_k\Symp_\tau(M,\pci) = 
\begin{cases}
H^{2-k}(M;\B Z) & \text {for $k\in \{1,2\}$}\\
               0 & \text {otherwise}
\end{cases}$
and hence
$\pi_k({\bukfa B}) \cong \pi_k(\BS\Mwo)$ for $k>3$. 
\end{corollary}

\begin{proof}
Since $\Map(M,\pci)$ is an H-space,
all its connected components are homotopy equivalent.
The first statement follows immediately from Theorem \ref{T:isomorphism}
and Lemma \ref{L:spanier}.
The second is the direct application of the long exact sequence of
homotopy group for the fibration 
$\Symp\Mwo\to\Symp_\tau(M,\pci)\to{\bukfa B}$.
\end{proof}

\begin{example}\label{E:pi0}
Let $\Si_g $ be a surface of genus $g>0$.
Let $i_0,i_1\in \Symp(\Si_g,\pci )$ be two embeddings such
that $i_0 = i_1\circ f$, where $f\colon \Si_g\to \Si_g$ is a
Dehn twist. One can construct a symplectic Lefschetz
fibration $\Mo \to S^2$ with the generic fiber $\Si_g$ such
that the monodromy about some critical value is exactly
$f$ and $\Mo $ is an integral symplectic manifold 
\cite{MR2002g:57051}. 
Embedding $\Mo$ symplectically into $\cp^N$ we also get a family
of symplectic embeddings $i_t\colon \Si_g\to \cp^N$ for $t\in [0,1]$.
Since $\pi_0(\Symp(\Si_g,area))$ is generated by Dehn
twists, we obtain that the map
$\pi_0(\Symp(\Si_g,area))\to \pi_0(\Symp(\Si_g,\pci))$
is trivial.
\end{example}

\begin{lemma}\label{L:pi2}
The map $\pi_2(\Symp\Mo) \to \pi_2(\Symp(M,\pci))$ is trivial.
In particular, 
$\pi_3({\bukfa B}) \cong \pi_3(\BS\Mo )$. 
\end{lemma}

\begin{proof}
Let $\xi \in \pi_2(\Symp(M,\om))$. Its image 
in $\pi_2(\Map(M,\pci))=\B Z$ is equal to 
$\left < [\om],ev_*(\xi)\right >$, where 
$ev\colon \Symp\Mo \to M$ is the evaluation at the base point $pt\in M$
(cf.~proof of Lemma \ref{L:spanier}). More precisely,
consider the following commutative diagram.
$$
\xymatrix
{
\Symp\Mo \ar[r]^i \ar[d]^{ev} & \Symp(M,\pci)\ar[d]^{ev}\\
M \ar[r]^f                    & \pci
}
$$
The isomorphism $\pi_2(\Map(M,\pci))=\B Z$ is given by
the integration of the generator $[\om_{\infty}]\in H^2(\pci)$
over the sphere $S^2 \to \Map(M,\pci)\stackrel {ev}\to \pci$.
Hence we calculate that
$
\left<[\om_{\infty}],ev_*i_*(\xi)\right>=
\left<[\om_{\infty}],f_*ev_*(\xi)\right>=
\left<f^*[\om_{\infty}],ev_*(\xi)\right>=
\left<[\om],ev_*(\xi)\right>
$
as required.

We claim that this number is equal to zero.
Indeed, if $ev^*[\om]$ is a sum of products of 1-dimensional classes
then $\left<ev^*[\om],\xi\right > = 
\langle \xi^*(ev^*[\om]),[S^2]\rangle = 0$. Suppose, by contradiction,
that $ev^*[\om]$ is not a sum of products.
Recall that, according to Milnor and Moore \cite{MR0174052},
the cohomology ring of $\Symp\Mo $ (which is an H-space) is
a free graded commutative algebra. If $ev^*[\om]$ were nonzero
then it would be of infinite order, which is impossible since
$[\om]^{n+1}=0$.
\end{proof}

Let $\Flux\colon \pi_1(\Symp\Mo)\to H^1(M;\B R)$ be the {\em flux homomorphism}
defined as follows. 
$\left <\Flux{\xi},\ga\right> = \left <[\om],\xi_*(\ga)\right>$,
where $\xi_*(\ga)\colon \B T^2\to M$ sends $(t,s)\to \xi_t(\ga(s))$ (see
Chapter 9 in \cite{MR2000g:53098} for details). Its image 
is called the flux group and is denoted by $\Ga_{\om}$.
Notice that in the case of an integral symplectic structure
$\Ga _{\om}$ is a subgroup of the integral cohomology $H^1(M,\B Z)$.
The following lemma is an immediate consequence  of the above definition.

\begin{lemma}\label{L:flux}
The following diagram is commutative.
$$
\CD
\pi_1(\Symp\Mo) @>>> \pi_1(\Symp(M,\pci))\\
@V\Flux{}VV           @V\cong VV\\
\Ga_{\om}       @>>> H^1(M,\B Z)=\pi_1(\Map(M,\pci))
\endCD
$$
\qed
\end{lemma}

\begin{proof}[Proof of Theorem \ref{T:BS}]
\hfill
\begin{enumerate}
\item
This part follows from Corollary \ref{C:vanish} and Lemma \ref{L:pi2}.
\item
To prove that the first sequence is exact
one applies Lemma \ref{L:pi2} and \ref{L:flux}
to the exact sequence associated with the fibration 
$\Symp\Mwo\to\Symp_\tau(M,\pci)\to{\bukfa B}$.
Similarly, Corollary \ref{C:vanish} and Lemma \ref{L:flux} applied
to the long exact sequence of the the same fibration proves the
exactness of the second sequence.
\end{enumerate}
\end{proof}

\subsection{Proof of Theorem \ref{T:gg}}\label{SS:gg}

Recall that principal $G$-bundles over $M$ are classified by
$H^1(M,G)$ and $H^1(M,U(1))$ is canonically isomorphic to $H^2(M,\B Z)$. 
The bundle corresponding to $\tau\in H^2(M,\B Z)$ will be denoted $L_\tau$,
and the class of a bundle $L$ in $H^2(M,\B Z)$ will be denoted by $[L]$.

Assume that $(M,\omega)$ and $(W,\omega_W)$
are integral symplectic manifolds.
Choose $\tau\in H^2(M,\B Z)$ and $\tau_W\in H^2(W,\B Z)$
to be some lifts of $\om$ and $\om_W$ respectively.
Let
$$\Map_{Symp}^{U(1)}(L_\tau,L_{\tau_W})$$
be the space of maps from $L_\tau$ to $L_{\tau_W}$
which are $U(1)$-equivariant and cover a symplectomorphism 
$\Mo\to \Wo$:
$$
\CD
L_\tau @>\hat\psi >> L_{\tau_W}\\
@VVV           @VVV\\
M @>\psi >>     W.
\endCD
$$
Observe that if $\psi\colon M\to W$ is a continuous map,
then $\psi^*(L_\tau)=L_{\psi^*\tau}$, thus there exists a covering
$U(1)$-equivariant map $\widehat f\colon L_{\psi^*\tau}\to L_\tau$.
In particular, for $\Wo =\Mo$ we get a group 
$$\C G:=\Map_{Symp}^{U(1)}(L_\tau,L_\tau).$$

\begin{lemma}\label{L:gg}
There is an  extension
$$
1\to \Map(M,U(1))\to \C G\to \Symp\Mwo \to 1
$$
\end{lemma}
\begin{proof}
The kernel of the map
that assigns to $\tilde\psi$
the underlying symplectomorphism $\psi$ consists of
the automorphisms
which cover the identity. Over each point of $M$ such
automorphism is given by an element of $U(1)$.
\end{proof}

\begin{proposition} 
Let $\om_n$ denote the standard symplectic form on $\cp^n$,
and let $[L_n]=[\om_n]\in H^2(\cp^n,\B Z)$.  
Then the space
$$\Map_{Symp}^{U(1)}(L_\tau,L_{\infty}):=
\bigcup_n \Map_{Symp}^{U(1)}(L_\tau,L_n)$$ 
is contractible.
\end{proposition}

\begin{proof}
Let $\Map^{U(1)}(L_\tau,L_\infty)$ be a space
of all $U(1)$-equivariant maps $L_\tau\to L_\infty$.
Then we have a bundle:
$$\Map(M,U(1))\to \Map^{U(1)}(L_\tau,L_\infty)
\to \Map_\tau(M,\pci),
$$
Consider the following two exact sequences of homotopy groups
$$
\xymatrix
{\dots \pi_k\left(\Map\left(M,U(1)\right)\right)\ar[r] &
\pi_k(\Map^{U(1)}(L_\tau,L_\infty))\ar[r] &
\pi_k(\Map_\tau(M,\pci)) \dots\\
\dots \pi_k\left(\Map\left(M,U(1)\right)\right)\ar[r]\ar[u]^=&
\pi_k(\Map^{U(1)}_{Symp}(L_\tau,L_\infty))\ar[r]\ar[u]&
\pi_k(\Symp_\tau(M,\pci))\ar[u] \dots
}$$
where the vertical arrows are induced by inclusions.
According to Theorem \ref{T:isomorphism}, the map
$$\pi_k(\Symp_\tau(M,\pci))\to\pi_k(\Map_\tau(M,\pci))$$
is an isomorphism. By the five lemma
the same is true for
$$\pi_k(\Map^{U(1)}_{Symp}(L_\tau,L_\infty))\to
\pi_k(\Map^{U(1)}(L_\tau,L_\infty)).$$

We are left to show that
$\Map^{U(1)}(L_\tau,L_\infty)$
is contractible. Since $L_\infty=EU(1)$,
the following lemma finishes the proof.
\end{proof}

\begin{lemma}
Let $E$ be a free $G$-space.
Then $\Map^G(E,EG)$ is contractible.
\end{lemma}

\begin{proof}
It is a standard fact (about universal principal
bundles) that $\Map^G(E,EG)$ is connected
(see Theorem 8.12 on page 58 in \cite{MR89c:57048}).
We need to show that for all $k$ the space
$\Map(S^k,\Map^G(E,EG))=\Map^G(E\times S^k,EG)$
is connected.
But  $S^k\times E$ is again a free $G$-space
and we apply the observation we began with.
\end{proof}

\begin{proof}[Proof of Theorem \ref{T:gg}]
Consider the following two fibrations
$$\Map(M,U(1))\to \Map_{Symp}^{U(1)}(L_\tau,L_\infty)
\stackrel{p_1}\to \Symp_\tau(M,\pci)$$

and
$$\Symp\Mwo \to \Symp_\tau(M,\pci)\stackrel{p_2}\to{\bukfa B}_\tau,$$
where the base ${\bukfa B}_\tau$ is defined to be this quotient.
Observe that the composition $p_2\circ p_1$ defines a fibration
$$
\CD
\C G\to  \Map_{Symp}^{U(1)}(L_\tau,L_\infty)
@>{p_2\circ p_1}>>
{\bukfa B}_\tau
\endCD
$$
with a contractible total space.
\end{proof}

\subsection{Proof of Theorem \ref{T:main} $(1\Rightarrow 2)$}
\label{SS:1=>2}

By assumption the connection form $\Om$ has a lift $\C T\in H^2(E,\B Z)
=[E,K(\B Z,2)]$. Since $K(\B Z,2) = \pci$, we 
get an embedding $f:E\to \cp^N$ such that $\C T$ is
the pull-back of the generator.
According to the original Gromov's h-principle (Theorem 3.4.2 A in
\cite{MR90a:58201}) we can deform $f$ to an immersion preserving forms.
That is an immersion $f'\colon E\to \cp^N$ so that $f^*(\om_N) = \Om$.
At the price of increasing $N$ we can deform it further
to get a embedding as we did in the proof of Theorem
\ref{T:isomorphism}.\qed

\subsection{Proof of Theorem \ref{T:main} $(1\Rightarrow 3)$}
\label{SS:1=>3}

A class $\C T$ defines a line bundle $L_{\C T}\to E$.
Notice that the pull-back bundle under the inclusion of the fiber
$i\colon M\to E$ is the line bundle $L_{i^*\C T}\to M$.
Thus $i^*\C T$ is a lift of $[\om]$ to $H^2(M,\B Z)$.

Consider the composition
$L_{\C T}\to E\to B$ as a bundle over $B$ with fiber $U(1)$.
The structure group of this fibration is contained in
$\C G$. Thus we get a classifying map
$\widehat f\colon B\to B\C G$
such that the following diagram commutes up
to homotopy (here $f$ classifies $E$):
$$
\xymatrix{
 &B\C G\ar[d]\\
B \ar[ur]^{\widehat f} \ar[r]^{f\phantom{dupaga}} & \BS\Mo
}
$$
\qed

Recall that, according to Theorem \ref{T:gg},
$B\C G\simeq \bukfa B$.
The classifying map $\widehat f\colon B\to B\C G={\bukfa B}$
is the one constructed in the previous section.

\section{The geometry of symplectic configurations}\label{S:gsc}

\subsection{A principal connection on $\Symp(M,W)$}\label{SS:connection}

The main reference for connections is Kobayashi-Nomizu 
\cite{MR97c:53001a}. The tangent space to $\Symp(M,W)$ at
point $f$ is equal to 
$$\{\C X\in \Ga (f^*TW)\colon\C X(p)=\frac {d}{dt} f_t(p) \quad f_t\in 
\Symp(M,W)\}.$$
We define the horizontal space $\C H_f$ to be the space of
the sections $f^*\om_W$-orthogonal to $M$:
$$\C H_f:=
\{\C X\in T_f\Symp(M,W)\colon f^*\om_W(\C X,Y)=0 \quad 
\forall Y\in \Ga(TM)\}.$$

Consider the following one--form $\theta $ on $\Symp(M,W)$
with values in closed one--forms on $M$ which is identified
with the Lie algebra of $\Symp\Mo$
$$
\theta_f (\C X):= f^*(\iota_{\C X_f}\om_{W})
$$
for $\C X\in T_f\Symp(M,W)$.
More precisely
$$
(\theta_f(\C X))(p)(Y) = \om_W(\C X_f(p),f_* Y(p)),
$$
where $Y$ is a vector field on $M$ and $\C X_f$ is a vector
field on $W$ defined on $f(M)$.

\begin{lemma}
The one--form $\theta$ is a connection form induced by $\C H$, 
that is, it satisfies
\begin{enumerate}
\item
$\theta_{f\circ \psi} = Ad_{\psi}\theta = \psi^*\circ \theta$,
\item
$\theta_f(\underline X) = X$, where $X\in \Om^1(M)$ is a closed one--form
and $\underline X$ denotes the fundamental vector field
(i.e. the vector field generated by the infinitesimal action of 
a Lie algebra),
\item
$\C H$ is the kernel of $\theta $.
\end{enumerate}
\end{lemma}

\begin{proof}
(1) This is an immediate calculation:
\begin{eqnarray*}
(\theta_{f\circ \psi}(\C X))(p)(Y) &=& 
\om_W\left(\C X_{f\circ \psi}(p), (f\circ \psi)_* Y(p)\right)\\
&=& 
\om_W\left(\C X_f (\psi(p)), f_*(\psi_* Y(p))\right)\\
&=&
(\theta_f(\C X_f))(\psi(p)) (\psi_*Y(p))\\
&=&
(\psi^*\circ \theta_f)(\C X_f)(p)(Y).
\end{eqnarray*}
(2) If $X$ is a closed one--form on $M$ then we denote by
$X^{\sharp}$ the $\om$-corresponding vector field.
That is $\iota_{X^{\sharp}}\om = X$.
Thus $\underline X_f = X^{\sharp}$ under the identification
between $M$ and $f(M)\in W$ given by $f$. Hence
we get
$$
\theta_f(X^{\sharp}) = f^*(\iota_X^{\sharp}\om_{W})= X.
$$

\NI
(3) Obvious.

\end{proof}

\begin{proposition}\label{P:sc}
The connection associated to $\theta$ on $M_W$ and connection
defined by the two--form $\Om_W:=ev^*(\om_W)$ coincide.
\end{proposition}

\begin{proof}
Both distributions consist of vectors $\C X\in T_{[f,p]}M_W$
satisfying the condition 
$ev_*(\C X)\in T_{f(p)}W$ is $\om_W$-orthogonal to $f_*(T_pM)$.
\end{proof}

\bs

\subsection{The curvature of $\theta$}\label{SS:curva}
Let's calculate the curvature two--form of the connection $\theta$.
By definition it is
$$
\Theta(\C X,\C Y) := d\theta(\C X^h,\C Y^h)= - \theta([\C X^h,\C Y^h]),
$$
where $\C X^h,\C Y^h\in \C H$
are horizontal parts of tangent vectors. Calculating it further we get
\begin{eqnarray*}
-\Theta_f(\C X^h,\C Y^h) &=& \theta_f[\C X^h,\C Y^h]\\
&=& f^*(\iota_{[\C X^h,\C Y^h]}\om_W)\\
&=& f^*(\{\al_{\C X},\al_{\C Y}\})
\end{eqnarray*}
Here $\al_{\C X},\,\al_{\C Y}$ are one--forms 
associated with extensions of $\C X^h,\,\C Y^h\in \Ga(f^*TW)$
to vector field on some neighbourhood of $f(M)\subset W$ and
$\{\al_{\C X},\al_{\C Y}\}$ is the Poisson bracket.
Since $\C X^h$ and $\C Y^h$ are horizontal, the 
one--forms $\al_{\C X},\,\al_{\C Y}$ vanish on $M$.
Therefore their local extensions might be chosen
to be exact and we get
$$
-\Theta_f(\C X^h,\C Y^h) = 
d\, \left [\,f^*\{H_{\C X},H_{\C Y}\}\,\right],
$$
where $H_{\C X},\,H_{\C Y}\colon W\to \B R$ are Hamiltonians,
i.e. $dH_{\C X}=\al_{\C X}$ and $dH_{\C Y}=\al_{\C Y}$.

In other words the curvature of the connection
$\theta$ at the embedding $f\colon \Mo \to \Wo$
is measured by the subspace
of functions on $M$ consisting of restrictions of Poisson 
brackets of functions on $W$ which are constant
on $f(M)$. We will use this curvature form to construct
characteristic classes in Section \ref{SS:Chern}.

\begin{lemma}\label{L:curv}
Let $f\colon \Mo\to\Mo \times (\Sigma ,\om_{\Si})$ 
be a symplectic embedding given by $f(p)=(p,s_0)$,
where $s_0\in \Si$ is a base point on the surface.
The image of the curvature two--form
$$\Theta_f\colon \C H_f\times \C H_f \to \text{Lie}\,(\Symp\Mo)$$
contains the subspace of all exact one--forms.
\end{lemma}

\begin{proof}
What in fact we need to show is that given a function
$H\colon M\to \B R$ there are functions 
$F,G\colon M\times \Si \to \B R$ constant
on $f(M)$ such that $\{F,G\}\circ f = H$.

Let $F',G'\colon \Si\to \B R$ be smooth functions 
such that $F'(s_0) = 0$ and 
$\{F',G'\}(s_0) = 1$. Define $F(p,s):=F'(s)$,
$G(p,s)=G'(s)$ and $\widehat H(p,s) = H(p)$.
The following calculation follows from the Leibniz property
of the Poisson bracket.
$$
\{\widehat HF,G\} = \widehat H\{F,G\} + F\{\widehat H,G\} = 
\widehat H\{F,G\}
$$
Hence, $\{\widehat HF,G\}\circ f = \widehat H\circ f = H$.
Note that the functions $\widehat H F$ and $G$ are constant
on $f(M)$ as required.
\end{proof}

\bs

\begin{corollary}\label{C:curv}
The image of the curvature two--form 
$$\Theta_{i\circ f}\colon 
\C H_{i\circ f}\times \C H_{i\circ f} \to \text{Lie}\,(\Symp\Mo)$$
of the connection in the universal configuration fibration
contains the subspace of all exact one--forms for every $f\colon \Mo\to \pci$. 
\qed
\end{corollary} 
\bs

\subsection{The holonomy group}\label{SS:holo}

Let $\Mo \stackrel i\to E \stackrel \pi\to B$ be a symplectic fibration.
Let $K:=\ker \,[\,\pi_*\colon H_2(E,\B Z)\to H_2(B,\B Z)]$.
The following exact sequence:
$$ 0\to H_2(M,\B Z)\to K\to H_1(B;H_1(F,\B Z))\to 0$$
is derived from the Leray-Serre spectral sequence associated
with $E$.
An element of $H_1(B,H_1(M,\B Z))$ may be described as
a (class of a) loop $\gamma$ in $B$ and a $\gamma$--invariant
cycle $[\ell]\in H_1(M,\B Z)^\gamma$ in the fiber 
$F$ over a base point of $\gamma$.

Let $\Om$ be a closed extension of $\om$. 
We construct a certain lift $\Sigma_{\gamma,\ell}\in H_2(E,\B Z)$ of 
$(\gamma,\ell)\in H_1(B;H_1(M,\B Z))$ to $K$. 
It is a union
of two chains $f(S^1\times I)$ and $C$. The map $f$ is such that
$f(\{s\}\times I)$ is a horizontal lift of $\gamma$ for each $s$ and
$f(S^1\times\{0\})=\ell$. The chain $C$ is a chain in $F$ such that
$\partial C=\gamma_*\ell-\ell$, where $\gamma_*$ is the holonomy of $\gamma$
with respect to $\Om$.
Observe that
$\Om(f_*(\frac d {ds}),f_*(\frac d {dt}))=0$ because the kernel
of $\Om$ is the horizontal subspace.
Thus $\int_{\Sigma_{\gamma,\ell}}\Om=\int_C\om$.

Let $\psi\in\Symp\Mo$ and
$H_1(M,\B Z)^\psi=\{\ell\in H_1(M,\B Z)\colon \psi_*(\ell)=\ell\}$
be the $\psi$--invariant part of the first cohomology.
Define a flux-like homomorphism
$$
\Flux^{\psi}\colon H^1(M,\B Z)^{\psi}\to \B R/\om
$$
as follows (here $\B R/\om$ denotes the quotient of $\B R$ by
the group of periods of the symplectic form).
If $[\ell]\in H_1(M,\B Z)^\psi$ then there exists a two--chain $C$
such that $\partial C=\ell-\psi_*\ell$. Then
$$\Flux^{\psi}(\ell):= \int_C\om\ \mod (\text{periods of } \om).$$
The integral is defined up to the period of $\om$
and does not depend on the choice of representative $\ell$ in the
homology class.

Let $\overline{\C D}$ be a subset of $\Symp\Mwo$ consisting of
those $\psi\in Symp\Mwo$ for which $\Flux^{\psi}=0$.
Note that $\overline{\C D}$ is conjugacy invariant,
but it is not a subgroup of $\Symp\Mwo$.

\begin{proposition}\label{P:flux}
Let $M\to E\to B$ be a symplectic fibration.
Let $\Om$ be an extension of the symplectic form $\om$ and
let $K:=\ker H_2(E,\B Z)\to H_2(B,\B Z)$.
Then  $\left <\Om, K\right> \subset \B Z$ if and only if
the holonomy $\pi_1(B)\to\Symp\Mwo$ is contained in $\overline{\C D}$.
\end{proposition}
\begin{proof}
If $C'\in H_2(M,\B Z)$ then $\int_{i_*(C')}\Om=\int_{C'}\om$
since $\Om$ extends $\om$.
If $(\gamma,\ell)\in H_1(B;H_1(M,\B Z))$ then
$$\int_{\Sigma_{\gamma,\ell}}\Om=\int_C\om=\Flux^{\psi}(\ell)
\ \mod (\text{periods of } \om),$$
where $\psi$ is the monodromy around the loop $\gamma$.
\end{proof}

Recall that we define $\C D$ to be the holonomy group of the
universal connection~$\theta$.
Notice that $\C D$ is not well defined as a subgroup of
$\Symp\Mo$ since it depends of the reference point
$f\colon \Mo \to \pci$. However, according to standard theory
of connections, different holonomy subgroups are conjugate.

\begin{corollary}\label{C:DD}
The set $\overline{\C D}$ is the sum of all subgroups conjugate to $\C D$.
\end{corollary}

\begin{proof}
Let $E=M\times\B R/\sim$, where $(m,t+k)\sim(\psi^k(m),t)$ for $k\in\B Z$.
The form $\Om=\pi^*_M\om$ on $M\times\B R$ descends
to a closed connection form on $E$ with a holonomy $\psi$.
Since $H_2(S^1,\B Z)=0$, we deduce from Proposition \ref{P:flux}
that if $\psi\in \overline{\C D}$ then
$\Om$ has integral periods.

To prove the converse inclusion notice that it is enough
to show that $\C D\subset \overline {\C D}$, according to the
conjugacy invariance.
Since $\C D$ is the universal holonomy group,
for any $\psi\in \C D$ there exists a bundle
$M\to E\to S^1$ with a closed integral connection form
$\Om$ such that $\psi$ is a holonomy of that form. Thus
$\psi\in\overline{\C D}$ again by Proposition \ref{P:flux}.
\end{proof}

\begin{proposition}\label{P:overline}
\hfill
\begin{enumerate}
\item
Let $H\subset\Symp\Mwo$ be any subgroup. If $H\subset\overline{\C D}$
then $H<\C D$ (up to conjugacy).
\item
$\C D\cap\Symp_0\Mwo=\overline{\C D}\cap\Symp_0\Mwo$,
\item
$\C D$ is a closed subgroup of $\Symp\Mwo$
\end{enumerate}
\end{proposition}
\begin{proof}
First we will prove (1) with the additional assumption that
$H$ is countable. The we will use it to prove (2) and (3).
Finally we will conclude general case of (1) with the use of (3).

\begin{enumerate}
\item
We shall construct a fibration over a noncompact
surface (possibly of infinite genus)
with closed connection form $\Om$ and holonomy equal to $H$. 
Since $H<\overline{\C D}$ then, according to Proposition \ref{P:flux}
$\left <\Om,K\right> \subset \B Z$. Because the base is
a noncompact surface then the latter is equivalent to 
the integrality of $\Om$. To see this notice that 
$K=H_2(E,\B Z)$ as $H_2(B,\B Z)$ is trivial.
By the universality of $\C D$ we will conclude that,
up to a conjugation, $H<\C D$.

Take $M_{\psi_i} \times S^1 \to T^2$ for each $\psi_i\in H$,
$i=1,2,\ldots$  and
form a countable fiber connected sum
$(M_{\psi_1} \times S^1) \# _M (M_{\psi_2} \times S^1) \# _M \dots.$
It is a bundle over $\Si_\infty$ the
surface of infinite genus and admits a closed connection form.
Indeed, every $M_{\psi}\times S^1$ admits
a closed connection form as in the proof of Corollary \ref{C:DD}.
One can easily extend it over the fiber connected sum.

Every loop in the connected sum is a composition of loops in the summands,
thus the holonomy of the bundle is contained in $H$. Since
all the generators of $H$ are contained in the holonomy, $H$ equals the
holonomy group of that bundle.

\item
Let $\psi \in \overline {\C D}\cap \Symp_0\Mo$. Observe
that $\Flux ^{\psi}\colon H^1(M,\B Z) \to \B R/\B Z$ is equal to
the value of the usual flux on $\psi$. This shows that
$\overline {\C D}\cap \Symp_0\Mo \subset \Flux^{-1}(H^1(M,\B Z))$.
In fact, the equality holds, because $\C D_0 = \Ham\Mo$,
according to Theorem \ref{T:gd} (2).

Now take one element $\psi$ from each component of
$\overline {\C D}\cap \Symp_0\Mo$ and generate subgroup
$H\subset \Symp_0\Mo$. It is clearly countable, hence
it is contained in $\C D\cap \Symp_0\Mo$, due to the 
part (1). Since $\C D_0 = \Ham\Mo$, we have
$\overline {\C D}\cap \Symp_0\Mo
=H\cdot \Ham\Mo \subset \C D\cap \Symp_0\Mo$.

\item
Since $\overline{\C D}\cap\Symp_0\Mwo = \Flux^{-1}(H^1(M,\B Z))$
is a closed subgroup of $\Symp_0\Mo$, (3) follows.
\end{enumerate}

To finish the proof we need to show that (1) holds in
full generality. Let $H'$ be a countable dense subgroup of $H$.
Such a group exists because $\Symp\Mo$ is a group modelled
on a separable and metrizable space of closed one--forms.
Thus it is second countable. The subgroup
$H$ is second countable as well, according to the hereditary 
properties of the second countability, and hence it is
separable. Take $H'$ to be a 
countable subgroup of $H$ generated by a countable 
dense subset.
Then $gH'g^{-1}\subset \C D$ by already proved part of (1)
and $gHg^{-1}\subset\C D$ by (3). 
\end{proof}

\subsection{Other results}\label{S:symbukfa}
We conclude this section with several geometric results.
The proofs are straightforward
and are left to the reader.

\begin{proposition}\label{P:haller}
\hfill
\begin{enumerate}
\item
The action of $\Symp\Wo$  on  
$\Symp(M,W)$ given by the composition
$\varphi \cdot f := \varphi \circ f$
preserves the connection $\theta$.
\item
The induced action on the associated bundle
$\Mo\to M_W\to {\bukfa B}_W$ is given by
$\varphi \cdot [f,p]:= [\varphi \circ f,p]$.
The evaluation
map $ev\colon M_W \to W$ is $\Symp\Wo$-equivariant.
\end{enumerate}
\end{proposition}

The general idea of this paper is to investigate
spaces ${\bukfa B}_W=\Symp(M,W)/\Symp\Mo$, 
where $\Mo$ is fixed and $\Wo$ varies. 
One can do the other way around and fix the target
manifold $\Wo$.
This approach was taken by Haller-Vizman in
\cite{haller}.
They write explicitly
a moment map for the action and express ${\bukfa B}_W$ as
a coadjoint orbit in the dual Lie algebra of the group
$\Ham\Wo$.

\begin{proposition}[\cite{haller}]\label{P:haller1}
\hfill
\begin{enumerate}
\item
The fiber integral $p_!(ev^*(\om_W^{n+1}))\in \Om^2({\bukfa B}_W)$
is a symplectic form on ${\bukfa B}_W$.
\item
The action induced on ${\bukfa B}_W$ 
preserves the symplectic form 
$$p_!(ev^*(\om_W^{n+1})).$$
\end{enumerate}
\end{proposition}

\section{Group cohomology and the group of McDuff}\label{S:group}

\subsection{Preparation on crossed homomorphisms}\label{SS:group}
Let $S$ be a topological group and $V$ be a topological $S$-module.
All the modules we consider are right modules.
A continuous map $\phi\colon S\to V$ is called a {\em continuous crossed
homomorphism} or {\em continuous 1-cocycle (with values in V)} if 
$\phi(gh)=\phi(h)+\phi(g)\cdot h$.

A crossed homomorphism $S\to V$ may be understood as a deformation
of the map $S\to V\rtimes_\rho S$. Precisely, the map
$$
(\psi,\rho)\colon S\to V\rtimes S
$$
is a homomorphism if and only if $\psi$ is a crossed homomorphism.

\subsection{An obstruction associated to a crossed homomorphism}\label{Bcrossed}

The kernel $H$ of a 1-cocycle $\psi$
is a pull-back in the following diagram:
$$
\CD
H@>>>S\\
@VVV @VVsV\\
S@>h>>V\rtimes S\\
\endCD
$$
Where $h=(\psi,\rho)$, and $s\colon S\to V\times S$ is the inclusion
onto the second factor.
Thus $BH$ is a homotopy pull-back in the following diagram
$$
\CD
BH@>>>BS\\
@VVV @VVBsV\\
BS@>Bh>>B(V\rtimes S)\\
\endCD
$$
The homotopy fiber of the map $Bs\colon BS\to B(V\rtimes S)$ equals  $V$,
moreover the fibration
$$
V\to BA=V\times_{V\rtimes S}E(V\rtimes S)\to B(V\rtimes S)
$$
is the universal $V$-bundle over $B(V\rtimes S)$.
Note that $V$ is an affine $V\rtimes S$-space.
In particular the obstruction to the existence of a section $BS\to BH$
is a pull-back of the obstruction to the existence of a section of
$Bs\colon BS\to B(V\rtimes S)$.

\begin{theorem}[McDuff]
There exists a class $\C O_M\in H^2(B\Symp\Mwo,H^1(M,\B Z))$ with the
following property. A $\Symp\Mwo$-bundle $M\to P\to B$ admits 
a reduction to $\C D_{\tau}$ 
if and only if $f^*\C O=0\in H^2(B,H^1(M,\B Z)/torsion)$,
where $f\colon B\to B\Symp\Mwo$ is the classifying map for the bundle $P$.
\end{theorem}

\begin{proof}
Let $S=\Symp\Mwo$ and $V=H^1(M,\B R/\B Z)$.
Due to Theorem \ref{T:dusa} (which we shall prove later),
$\C D_{\tau}$ is the kernel of the cocycle of McDuff 
$F^s\colon \Symp\Mwo\to H^1(M,\B R/\B Z)$. 
Since $V=K(H^1(M,\B Z),1)$, the total obstruction class
$\C O_V$ to the existence of the lift $B(V\rtimes S)\to BS$ is 
the element of $H^2(BS,\pi_1(V))=H^2(BS,H^1(M,\B Z))$.
Thus we set $\C O_M=(Bh)^*\C O_V$ where $h=(\rho,F^s)$ 
and  $\rho$ is an action of
$\Symp\Mwo$ on $H^1(M,\B Z)$.
\end{proof}

In Section \ref{SS:bockstein} we describe the Bockstein map $\beta$
which, for given short exact sequence of $S$-modules
$$
0\to Z\to R\to R/Z\to 0,
$$
to every (cohomology class of a) crossed homomorphism $S\to R/Z$ associates
an extension $0\to Z\to G\to S\to 0$.
If $Z$ is discrete, then such an extension is classified
by an element $\sigma_G\in H^2(BS,Z)$ (cf. \cite{MR43:6292}).

\begin{corollary}
$\C O_M=\sigma_{\C G/\Map_0(M,U(1))}=\beta([F^s])$.
\end{corollary}

\begin{proof}
By the result of \ref{SS:bockstein},
$\C D_{\tau}$ and $\C G/\Map_0(M,U(1))$ are homotopy
equivariant, thus the obstruction to the existence of a section
$\C D_{\tau}\to\Symp\Mwo$ equals the the obstruction to the existence
of a section $\C G/\Map_0(M,U(1))\to\Symp\Mwo$.

The second equality follows directly from the definition of $\beta$.
\end{proof}

\subsection{The homomorphism of Bockstein}\label{SS:bockstein}
Let $S$ be a topological group,
$R$ be a topological $S$-module, and $Z$ its submodule.
Let $\phi\colon S\to R/Z$ be a crossed homomorphism.
Define $G$ to be the subgroup of the semidirect product
$R\rtimes S$ consisting of the pairs $(s,r)$ such that $\phi(s)=r+Z$.
We call $G$ a pull-back of the following diagram:
$$
\CD
G @>\Phi>> R\\
@VpVV  @VVV\\
S @>\phi>> R/Z\\
\endCD
$$
where $p$ and $\Phi$ are the projections on the factors.
The vertical arrows are homomorphisms
and $\Phi$ is a crossed homomorphism with respect to the action
factored through $S$.

\begin{remark}
$G$ is a pull-back in the category theory sense \cite[p. 81]{MR1950475} of the
following diagram of groups
$$
\CD
G @>(\Phi,\rho\circ p)>>R\rtimes S\\
@VpVV  @VVV\\
S @>(\phi,\rho)>>(R/Z)\rtimes S\\
\endCD
$$
\end{remark}

\begin{lemma}\label{L:section}
Let $H$ be a kernel of $\phi$. Then there exist a section $\sigma\colon H\to G$
(i.e. continuous homomorphism, such that $p\circ\sigma=id_H$).
\end{lemma}

\begin{proof}
Consider a pull-back:
$$
\CD
G'@>\Phi'>> R\\
@Vp_HVV  @VVV\\
H @>\phi_{|H}>> R/Z\\
\endCD
$$
One observes that since $\phi_{|H}$ is trivial then $G'$ is a semidirect
product $Z\rtimes H$ and, in particular, there is a splitting $H\to G'$.
Clearly, $G'$ is a subgroup of $G$. Hence composing the splitting
$H\to G'$ with the inclusion $G'<G$ we get the required section.
\end{proof}

Note that $\ker p = Z$ and consider the two (topological) fibrations:
$$ Z\to G\stackrel p\to S,$$
and 
$$H\to S\stackrel\phi\to R/Z.$$

\begin{lemma}\label{L:longsequence}
The corresponding long exact sequences of homotopy groups are
matched up as in the following commutative diagram:
$$
\CD
\cdots@>>>\pi_{k+1}(S)@>>>\pi_{k+1}(R/Z)@>>>\pi_k(H)@>>>\cdots\\
&&@| @V\partial_RVV @V\pi_k(\sigma)VV &&\\
\cdots@>>> \pi_{k+1}(S)@>>>\pi_k(Z)@>>> \pi_k(G)@>>>\cdots\\
\endCD
$$
where $\partial_R$ is the connecting homomorphism
in the long exact sequence of the homotopy groups of the
fibration $Z\to R\to R/Z$.
\end{lemma}

\begin{proof}
This is straightforward to check using standard diagram-chasing. 
\end{proof}

\begin{corollary}\label{P:group_he}
If $R$ is contractible then $\sigma$ constructed above
is a homotopy equivalence.
\end{corollary}

\begin{proof}
If $R$ is contractible then $\partial_R$ is an isomorphism. 
Hence, by the five lemma, $\pi_k(\sigma)$ is also an isomorphism
and the statement follows.
\end{proof}


\subsection{The cocycles and group of McDuff}\label{SS:cocycle_dusa}
Let $\om$ be a closed (possibly singular) two--form on a manifold $W$.
Let $SH_1(W;\B Z)$ be space of integral 1-cycles
quotient out by boundaries of chains of zero $\om$-area. 
There is a split extension
$$
0\to\B R/\om \to SH_1(W,\B Z)\stackrel p \to H_1(W,\B Z)\to 0.
$$
The existence of a section $s\colon H_1(W,\B Z)\to SH_1(W,\B Z)$
follows since $\B R/\om$ is divisible. 
Precisely such a section may be constructed inductively, as follows.
Assume that the section is defined on a submodule $V\subset H_1(W,\B
Z)$,
$v\in H_1(W,\B Z)$ and $n$ is the least integer such that $nv\in V$.
Let $\tilde v$ by any lift of $v$ in $SH_1(W,\B Z)$.
If $n=\infty$ we set $s(v)=\tilde v$. In the other case
we see that $s(nv)-n\tilde v=[r]\in\B R/\om$. If $\gamma_{r/ n}$
is a contractible loop bounding a disk of area $r/n$ we set
$s(v)=\tilde v +\gamma_{r/ n}$.

Clearly, any group $G$ that fixes $\om$ acts on $SH_1(W,\B Z)$.
In \cite{dusa} McDuff defines a continuous
crossed homomorphism $F_s\colon G\to H^1(M;\B R/\om)$ as follows
$$
F_s(g)(\ga):= s(g_*(\ga))-g_*(s(\ga))\in
\B R/\om=\ker\left(SH_1(W,\B Z)\to H_1(W,\B Z)\right).$$
\begin{lemma}
$F_s$ satisfies
$$F_s(gh)(\ga) = g_*(F_s(h)(\ga)) + F_s(g)(h_*(\ga)) .$$
\end{lemma}
\begin{proof} Write for short $F_s(g)=[s,g_*]$. Then
$$F_s(gh)=[s,g_*h_*]=g_*[s,h_*]+[s,g_*]h_*=g_*F_s(h)+F_s(g)h_*$$
by the Leibniz rule.
\end{proof}

Since $G$ acts trivially on $\B R/\om$, $g_*(F_s(h)(\ga))=F_s(h)(\ga)$.
Hence $$F_s(gh) = F_s(g) + (F_s(g))\cdot h,$$
that is $F_s$ is indeed a crossed-homomorphism.

We shall apply this construction in the following situation:
\begin{enumerate}
\item $(W_1,\om_1)=(L_\tau,\pi^*\om)$ where $L_\tau$ is a prequantum
  $U(1)$-bundle
with $[L]=\tau\in H^2(M,\B Z)$ and $G_1=\C G$.
\item $(W_2,\om_2)=\Mo$ and $G_2=\Symp\Mwo$.
\end{enumerate}

\begin{lemma}\label{L:cocycle}\hfill
\begin{enumerate}
\item Since $\pi^*\om$ is exact on $L_\tau$, the group of periods
  is trivial
i.e. $H^1(L_\tau,\B R/\pi^*\om)=H^1(M,\B R)$.
\item One can choose sections
$\tilde s\colon H_1(L_\tau,\B Z)\to SH_1(L_\tau,\B Z)$
and $s\colon H_1(M,\B Z)\to SH_1(M,\B Z)$
such that the following diagram commutes.
$$
\xymatrix{\B R\ar[r]\ar[d] & SH_1(L_\tau,\B Z)\ar[r]\ar[d] &
  H_1(L_\tau,\B Z)\ar[d]^=\ar@<-4pt>[l]_{\tilde s}\\
\B R/\om\ar[r]   & SH_1(M,\B Z)\ar[r] & H_1(M,\B Z)\ar@<-4pt>[l]_s}
$$
\item If $\phi\colon \C G\to \Symp\Mwo$ assigns to a bundle map
its underlying symplectomorphism then the following diagram commutes:
$$
\CD
\C G @>F_{\tilde s}>> H^1(M;\B R)\\
@V\phi VV                  @VVV       \\
\Symp\Mwo @>F_s>>          H^1(M;\B R/\om)
\endCD
$$
\item The restriction of $F_{\tilde s}$ to $\Map_0(M,U(1))$ is
  trivial,
and thus descends to a crossed homomorphism
$\widehat F_s\colon \C G/\Map_0(M,U(1))\to H^1(M;\B R)$.
\item The following diagram
$$
\CD
\C G/\Map_0(M,U(1))@>\widehat F_{\tilde s}>> H^1(M;\B R)\\
@VVV                  @VVV       \\
\Symp\Mo @>{F_s}>> H^1(M;\B R/\B Z)
\endCD
$$
is a pull-back.
\end{enumerate}
\end{lemma}

\begin{proof}
Most of the statements are trivial. 

\NI
(2) First we have to define $\tilde s$ on $K=\ker H_1(L_\tau,\B Z)\to
    H_1(M,\B Z)$
in such a way that its image lives in 
$\ker SH_1(L_\tau,\B Z)\to SH_1(M,\B Z)$.
The generator of $K$ is represented by a loop in a fiber, and
the value of $\tilde s$ may also be taken to be a loop in the
    fiber.
Then we extend $\tilde s$ to whole $H_1(L_\tau,\B Z)$ as explained
    above.
It clearly descends to the map $s\colon H_1(M,\B Z)\to SH_1(M,\B Z)$.

\NI
(4) This is obvious since $Map_0(M,U(1))$ is connected and the fiber
of the map $H^1(M,\B R)\to H^1(M,\B R/\B Z)$ is discrete.
\end{proof}

\begin{remark}[about Lemma \ref{L:cocycle} (2)]
McDuff \cite{dusa} shows that for each $s$ one can choose a lift
$\tau$ of $[\om]$ and a section $\tilde s$ such that $s$ is
constructed
as above.
\end{remark}

McDuff proves that the restriction of $F_s$
to $\Symp_0\Mo$ is equal to the flux homomorphism.
She defines $\Ham^{s\B Z}\Mo := \ker{F_s}$.
Clearly, the component of the identity of $\Ham^{s\B Z}\Mo$ 
is equal to $\Ham\Mo$.

\begin{corollary}\label{C:section}
There exists a continuous group homomorphism 
$$\Ham^{s\B Z}\Mo \to \C G/\Map_0(M,U(1))$$
that is a homotopy equivalence.
\end{corollary}

\begin{proof}
Since $H^1(M;\B R)$ is contractible, the statement follows
from Corollary \ref{P:group_he}.
\end{proof}

\subsection{The groups $\Ham^{s\B Z}$ and $\C D$ are equal}
Let $G$ be a topological group and 
$F\colon G_0\to A$ a $G$-equivariant
homomorphism to a right $G$-module $A$. That is it satisfies
$F(f^{-1}\, g \, f)= F(g)\cdot f$.
Let $H$ be a subgroup of $G$ intersecting every
connected component of $G$ and such that $H\cap G_0=\ker F$. 
Let $\sigma\colon G/G_0=\pi_0(G)\to  H$ be a set theoretic section of
the natural projection $p\colon G\to G/G_0$.
Define 1-cocycle $F^\sigma\colon G\to A$ by 
$$F^\sigma(g):= F(\sigma(p(g))^{-1}g).$$
The following lemma that is a reformulation of Proposition 1.8(ii) from
\cite{dusa}.

\begin{lemma}\label{L:subgroup}
The 1-cocycle $F^\sigma$ extends $F$ and
$H=\ker F^\sigma$.
\end{lemma}

\begin{proof}
It is straightforward that $F^\sigma$ extends $F$.
To prove that $F^\sigma$ is a 1-cocycle choose $f,g\in G$.
Since $\sigma(G)\subset H$,
$h=\delta(\sigma\circ p)(f,g)=
\sigma(p(fg))\sigma(p(g))^{-1}\sigma(p(f))^{-1}\in H$.
We calculate
\begin{eqnarray*}
G_0\ni\sigma(fg)^{-1}fg&=&h\sigma(p(g))^{-1}\sigma(p(f))^{-1}fg\\
&=&h\sigma(p(g))^{-1}gg^{-1}\sigma(p(f))^{-1}fg\\
&=&h\left(\sigma(p(g))^{-1}g\right)\left(g^{-1}\sigma(p(f))^{-1}fg\right)
\in hG_0.
\end{eqnarray*}
Thus $h\in G_0$ and finally $h\in\ker F=G_0\cap H$.
Then
\begin{eqnarray*}
F^\sigma(fg)&=&F\left(\sigma(p(f\,g))^{-1}fg\right)\\
&=&F(h)+F\left(\sigma(p(g))^{-1}g\right)+F\left(g^{-1}\sigma(p(f))^{-1}fg\right)
\\
&=& F^\sigma(g)+F^\sigma(f)\cdot g
\end{eqnarray*}

The proof that $H=\ker F^\sigma$ is also immediate and
is left to the reader.
\end{proof}

\begin{theorem}\label{T:dusa}\hfill
The group $\Ham^{s\B Z}\Mo$ and the holonomy group $\C D$ are
isomorphic.
\end{theorem}

\begin{proof}
{\it Step 1.} We will show that
$\Ham^{s\B Z}\Mo\subset\overline{\C D}$.
We have to observe that if $[l]\in H_1(M,\B Z)^\psi$
then $F_s(\psi)(l)=\Flux^\psi(\ell)$.

{\it Step 2.} By Proposition \ref{P:overline} (1) we conclude that
$\Ham^{s\B Z}\Mo\subset\C D$ (up to conjugacy).

{\it Step 3.} Let $F$ be the following composition
$$
\Symp_0\Mwo\stackrel{\Flux}\to H^1(M,\B R)/\Gamma_\om\to
H^1(M,\B R/\B Z).
$$
It follows from the definition that
$\Ham^{s\B Z}\Mo\cap\Symp_0\Mwo=\ker F$.
Let $\sigma$ be a section $\pi_0(\Symp\Mwo)\to\Ham^{s\B Z}\Mo$.
Then by Lemma \ref{L:subgroup} $\Ham^{s\B Z}\Mo=\ker F^\sigma$.
By Theorem \ref{T:gd}(3) $\C D=\Ham\Mo\cap\Symp_0\Mwo=\ker F$.
Since $\sigma(\pi_0(\Symp\Mwo))\subset \Ham^{s\B Z}\Mo\subset\C D$,
by Lemma \ref{L:subgroup}, we conclude that $\C D=\ker F^\sigma$.
\end{proof}

\subsection{Proof of  Theorem \ref{T:ugly}}
\label{SS:4=>3}

Let $\tau\in H^2(M,\B Z)$ be a preimage of the class 
$[\om]\in H^2(M,\B R)$ of the symplectic form. 
Let $\Mo \to E \to B$ be a symplectic fibration with the
classifying map $f\colon B\to B\Symp\Mwo$. 

\begin{lemma}\label{L:d_2}
There exists a lift
$\widehat f\colon B\to B\left(\C G/\Map_0(M,U(1))\right)$
if and only if 
$d_2(\tau) = 0$. 
Here $d_2$ is the differential in the spectral sequence associated
to $E$.
\end{lemma}
\begin{proof}


There exists a
lift $\widehat f\colon B\to B\left(\C G/\Map_0(M,U(1))\right)$
if and only there exist a lift $\widetilde f\colon B_2\to B\C G$
over the two--skeleton of $B$. 
Indeed,  we take a lift over the two skeleton by composing
the lift $\widetilde f$ with the projection
$B\C G\to B(\C G/\Map_0(M,U(1)))$.
Since the fiber $BH^1(M,\B Z)$
of the bundle $B(\C G/\Map_0(M,U(1)))\to B\Symp\Mwo$
has trivial higher homotopy
then all the further obstructions vanish and we extend
the lift to whole $B$.

Next, according to Theorem \ref{T:main} (2)$\iff $(3),
the existence of a lift $\widehat f\colon B_2\to B\C G$
is equivalent to the fact that the class 
$\tau=i^*(\widehat\tau)$, where $i\colon M\to E_2$ is the inclusion
of the fiber into the total space of the restriction of
$E\to B$ over the two--skeleton $B_2$. This in turn is equivalent
to the vanishing of $d_2(\tau)$ in the spectral sequence
for $E$.



\end{proof}

\begin{lemma}[\cite{MR2115670}]\label{L:d_3}
If $d_2(\tau) = 0 $ then $d_3(\tau)\in E_3^{3,0}$ is a torsion element.
\end{lemma}
\begin{proof}
Recall that $\dim M=2n$
and that $\pi_1(\BS\Mwo)$ acts trivially on $H^{2n}(M,\B R)$.
Then we have that
\begin{eqnarray*}
0&=&d_3[\tau^{n+1}]\\
&=& (n+1)d_3[\tau]\otimes[\tau^n]\in
H^3(\BS\Mo;H^{2n}(M,\B Z)).
\end{eqnarray*}
It means that $d_3[\tau]$ is a torsion element.
\end{proof}

\begin{proof}[Proof of Theorem \ref{T:ugly}]
Suppose that there is a lift $\widehat f\colon B\to B\C D$.
It implies, due to Theorem \ref{T:dusa}, that
there is a lift $\widehat f\colon B\to B\C G/\Map_0(M,U(1))$.
Lemma \ref{L:d_2} and Lemma \ref{L:d_3} show 
that $d_2(\tau)=0$ and $d_3(\tau)$ are 
torsion elements. 

Since all differentials of $[\om]$ vanish, 
there exists a closed extension $\Om'$ of $\om$.
We shall define a new closed connection form
$\Om:=\Om' + \pi^*\be$ so that it is integral (cf.
Section 3 in \cite{dusa}). Since $H_2(B,\B Z)$ is
torsion-free,
$H_2(E, \B Z) = \pi_*(H_2(E,\B Z)) \oplus K$, where 
$K:= \ker \pi_*$ as in Section \ref{SS:holo}.
It follows from Proposition \ref{P:flux}
that $\Om'$ has integral periods on the cycles from $K$.
We then define $\be\in \Om^2(B)$ to
be any two--form satisfying
$\int_{\pi_*(C_i)}\beta=-\int_{C_i}\Om'$,
where $\{\pi_*(C_i)\}$ forms a basis of
the image $\pi_*(H_2(E,\B R))$.

Finally, since $E$ admits closed integral connection form
it is an integral symplectic configuration, due to
Theorem \ref{T:main}. Notice that the class of the
above integral connection form might not restrict to
$\tau$, i.e. $\tau \neq i^*[\Om]\in H^2(M,\B Z)$.

\end{proof}

\section{The cohomology ring of symplectic configurations}\label{S:cohomology}

\subsection{Characteristic classes of configurations}

We define a characteristic class of configurations
of $\Mo$ in $\Wo$ to be an element of the cohomology ring
$H^*({\bukfa B}_W)$. If $M$ is compact then certain characteristic
classes can be obtained as fiber integrals 
$$\chi^W_k:=p_!(ev^*([\om_W^{n+k}])).$$
These classes are natural in the sense that if
$\psi\colon \Wo \to (V,\om_V)$ is a symplectic embedding then
$
\chi^W_k = \Psi ^*\chi^V_k,
$
where $\Psi\colon {\bukfa B}_W\to {\bukfa B}_V$ is the map induced
by $\psi$. 

Since any integral symplectic manifold embeds into $\cp^N$,
we define universal characteristic classes with respect
to $\pci$. We denote
the above fiber integrals  by $\chi_k\in H^{2k}({\bukfa B})$.
In this case those are usual characteristic classes of a fibration
with the structure group $\C G$.

A fundamental question is whether these classes are
{\em symplectic} that is whether they come from $H^*(B\Symp)$
via the canonical map ${\bukfa B}\to B\Symp$.
In general the answer is negative.

\begin{lemma}
Let $\si \in \pi_2({\bukfa B})$ be an element coming from
$H^0(M;\B Z)$ (see Theorem \ref{T:BS}). Then
$\left<\chi_1,\si\right>\neq 0$. In particular,
$\chi_1$ is not symplectic class. Moreover, $\chi_1$ is
of infinite order in $H^*({\bukfa B})$.
\end{lemma}

\begin{proof}
Recall that $\si$ is represented by a configuration
$M\times S^2 \to \pci$ such that the $S^2$ summand
is represents the positive generator of $\pi_2(\pci)$.
Clearly, $ev^*(\om_{\infty}) = \om \oplus \om_{S^2}$
and hence $\pi_!(ev^*(\om_{\infty}^{n+1}))=
(n+1)\,\,\om_{S^2}$.
Since the configuration is a trivial fibration,
no symplectic characteristic class can be nontrivial.

In order to prove the last statement, consider a configuration
$P:=M\times (S^2\times \dots \times S^2)\to \pci$ such that
every $S^2$ is as above.
Clearly, $\chi_1(P)$ is a sum of generators
and its top power is nonzero. Since it works for any number
of factors $S^2$, it follows that $\chi_1$ is of
infinite order.
\end{proof}

Let $\Mo\stackrel i\to M_{\Ham\Mo}\stackrel p\to \BH\Mo$ 
be the universal Hamiltonian fibration.
Recall that the {\em coupling class}
$\Om \in H^2(M_{\Ham\Mo})$ 
is defined by the following two conditions:
\begin{enumerate}
\item
$i^*(\Om) = [\om]$
\item
$p_!(\Om^{n+1})=0$
\end{enumerate}
One can define characteristic classes
by $\mu_k:=p_!(\Om^{n+k})\in H^{2k}(\BH\Mo)$ 
\cite{JK,MR2115670}.

Let $\widetilde{\bukfa B} := \Symp(M,\pci) / \Ham\Mo$ and
let 
$\Mo \stackrel i\to M_{\widetilde{\bukfa B}} 
\stackrel p\to \widetilde {\bukfa B}$
be the associated bundle. Notice that 
$\widetilde {\bukfa B}\to {\bukfa B}$ is the universal
cover. We have the following commutative
diagram of classifying spaces:
$$
\xymatrix{
\widetilde {\bukfa B} \ar[r]\ar[d] & {\bukfa B}\ar[d]\\
\BH\Mo \ar[r]                    & \BS\Mo
}
$$
We want to compare the pull-backs of characteristic classes
$\chi_k$ and $\mu_k$ in $H^*(\widetilde {\bukfa B})$. To avoid
clumsy notation we will denote these pull-backs by the same
symbols as the original classes.
The guiding idea is that the $\chi_k$ classes might
be easier to calculate than $\mu_k$'s.

\begin{lemma}\label{L:chi_k}
Let 
$\Mo \stackrel i\to M_{\widetilde{\bukfa B}} 
\stackrel p\to \widetilde {\bukfa B}$ be the above associated bundle.
Then $\chi_k \cong \mu_k $ modulo an ideal in $H^*(\widetilde {\bukfa B})$
generated by $\chi_1$.
\end{lemma}

\begin{proof}
First observe that the coupling class 
$\Om = ev^*(\om_0) - \frac{1}{n+1} p^*(\chi_1)$.
Calculating the appropriate fiber integral we get
the statement.

\begin{eqnarray*}
\chi_k&=&p_!(ev^*(\om_0^{n+k}))\\
&=& p_!\left((\Om + \frac{1}{n+1}p^*(\chi_1))^{n+k}\right)\\
&=&p_!\left(\sum_{i=0}^{n+k} 
\left(\begin{array}{c}n+k\\i\end{array}\right)
\Om^{n+k-i}\left(\frac 1 {n+1} p^*(\chi_1)\right)^i\,\,\,\right)\\
&=&\sum_{i=0}^{n+k} 
\left(\begin{array}{c}n+k\\i\end{array}\right)
\left(\frac {1}{n+1}\right)^i p_!(\Om^{n+k-i})\,\,(\chi_1)^i\\
&=&\sum_{i=0}^k 
\left(\begin{array}{c}n+k\\i\end{array}\right)
\left(\frac {1}{n+1}\right)^i \mu_{k-i}\,\,(\chi_1)^i
\end{eqnarray*}

\end{proof}

\begin{corollary}\label{C:chi_k}
Let $\Mo \to P\to B$ be a configuration whose classifying
map admits a lift to $\widetilde {\bukfa B}$.
If $H^2(B)=0$ then
$\chi_k(P) = \mu_k(P)$. In particular,
if $\chi_k(P)$ is nonzero then so does $\mu_k(P)$.
\qed
\end{corollary}

\subsection{Calculating $\chi$-classes}

It is important to know if there are relations between
characteristic classes. If $\Mo \to P \to S^{2k}$
is a symplectic configuration with nontrivial $\chi_k(P)$
then $\chi_k$ is not a product in $H^*({\bukfa B})$.
Indeed, if it were a product then its pull-back
$\chi_k(P)\in H^{2k}(S^{2k})$  would vanish.

\begin{theorem}\label{T:cha}
Let $\Mo \to P \to S^{2k}$ be a symplectic configuration, $k>1$.
The following conditions are equivalent.
\begin{enumerate}
\item
$\chi_k(P)\neq 0$;
\item
$\mu_k(P)\neq 0$;
\item
$ev_*[P]\neq 0$ in $H_{2n + 2k}(\pci;\B Z)$;
\item
There exists $a\in H^2(P)$ such that $a^{n+k}\neq 0$
(i.e.~$P$ is c-symplectic), and $a$ extends~$[\omega]$.
\end{enumerate}
Moreover, any of these conditions implies that $k\leq n+1$.
\end{theorem}

\begin{proof}
The first two conditions are equivalent, according to Corollary
\ref{C:chi_k}. Each of them is equivalent to the third
by the basic properties of fiber integration:
$$\left<\chi_k(P),[S^{2k}]\right> =
\left<ev^*(\om_{\infty}^{n+k}),[P]\right> =
\left<\om_{\infty}^{n+k},ev_*[P]\right>.$$
The third condition implies the fourth one since
the c-symplectic class is given by $[ev^*(\om_{\infty})]$.
The converse implication follows from
Theorem \ref{T:main}.

For the last statement, observe that $H^{2n+2}(P)=0$
if $k>n+1$. Hence $P$ cannot be c-symplectic.
\end{proof}

\begin{remark}
Consider the question of finding a large $k\in \B Z$
such that the Hurewicz map 
$\pi_{2k}(\BS\Mo)\to H_{2k}(\BS\Mo;\B Q)$
is nontrivial. The characteristic classes
constructed as fiber integrals give quite restricted
method for answering this question. As mentioned in the 
above theorem one can detect homologically nontrivial
spheres up to dimension $\dim M +2$.
\end{remark}

\subsection{Examples}
\hfill
\begin{enumerate}
\item
Consider the symplectic fibration 
$\cp^1 \to \cp^3\to S^4=\B H\B P^1$.
The total space is symplectic hence, due to Theorem \ref{T:cha},
we get that $\mu_2 = \chi_2 \neq 0$ in $H^4({\bukfa B}(S^2))$.

\item
Moreover, take a symplectic fibration
$F \to P \to \cp^3$ such that $P$ admits a compatible
symplectic structure. Composing with the above 
we get another symplectic fibration
$$
\Mo \to P \to S^4,
$$
whose fiber is the total space of 
$F \to \Mo \to \cp^1$. Again, we obtain that
$\chi_2 = \mu_2 \neq 0$ in $H^4({\bukfa B}\Mo)$.

\item
Similarly, consider the symplectic fibration
$$\cp^3 \to \cp^7 \to S^8=\B C{\text a}\B P.$$
Here, we get that $\chi_4$ and $\mu_4$ are nonzero
for $\cp^3$.

\item
Taking a symplectic fibration over $\cp^7$ admitting a
compatible symplectic structure we get a result
analogous to the one in (2).

\item
Let 
$\Mo :=
U(m) / U(m_1) \times \dots \times U(m_l), \,\,m=m_1+\dots +m_l$
and $m_1\geq m_2 \geq \dots \geq m_l$ 
be a generalized flag manifold equipped
with the homogeneous symplectic form. 
It is proved in K\c edra-McDuff \cite{MR2115670} that the classes
$\mu_k\in H^{2k}(\BS_0\Mo)$ are nonzero for $1<k\leq m$.
Moreover, it is shown that they are nontrivial for
certain fibrations over spheres, say
$\Mo\to P\to S^{2k}$, with the structure group
$SU(m)$. According to Theorem \ref{T:cha}, $P$ is c-symplectic,
however, it is not clear if it admits any symplectic form
in general. 

As in the previous examples we consider a symplectic fibration
$F\to E\to P$ and compose it with $P\to S^{2k}$. Suppose
that this composition is a symplectic fibration with fiber
$F\to Q\to M$. Then the classes $\mu_k(E)\in H^{2k}(\BS_0(Q))$
are nonzero.

\end{enumerate}

\subsection{Chern-Weil approach}\label{SS:Chern}

Recall from Section \ref{SS:curva} 
that the universal configuration fibration admits
a connection with curvature form
$\Theta$. It is an $ad$--invariant two--form on the total space
$\Symp(M,\pci)$ of the universal principal fibration
with values in the Lie algebra of $\Symp\Mo$ which
is identified with closed one--forms on $M$.
Since the reduced holonomy group
is Hamiltonian (Theorem \ref{T:gd} (1)), the curvature form has values in
the Lie algebra of $\Ham\Mo$ which is identified with
$C^{\infty}_0(M)$,
the space of functions on $M$ with zero mean value.

Let $\C P\colon C^{\infty}_0(M)\to \B R$ be an invariant polynomial.
It means that $\C P(h_1\circ \psi,...,h_k\circ \psi) = \C P(h_1,...,h_k)$,
where $\psi\in \Symp\Mo$ and $h_i\in C^{\infty}_0(M)$.
Then 
$$
\chi^{\C P}_{[f]}(X_1,Y_1,\dots ,X_k,Y_k):= 
\C P(\Theta_f(X_1^h,Y_1^h),\dots ,\Theta_f(X_k^h,Y_k^h))
$$
defines a $2k$-form on the base $\bukfa B$. The definition does not
depend on the choice of $f$ in the fiber since the polynomial
$\C P$ is invariant. It is a standard fact that these forms
are closed \cite{MR97c:53001b}. 
We define universal characteristic
classes to be cohomology classes of the form $\chi^{\C P}$.

The following formula gives a sequence of invariant polynomials.
$$
\C P_k(h_1,\dots ,h_k) := \int_M h_1\dots h_k\, \om^n
$$
The proof of the next lemma is analogous to the proof of
Lemma 3.9 in K\c edra-McDuff \cite{MR2115670}.
\begin{lemma}
$\left [\chi^{\C P_k}\right ] = \text{\em const } \chi_k$ for $k>1$.\qed
\end{lemma}

\bibliography{../../bib/bibliography}

\def\polhk#1{\setbox0=\hbox{#1}{\ooalign{\hidewidth
  \lower1.5ex\hbox{`}\hidewidth\crcr\unhbox0}}}
\begin{thebibliography}{ABKP00}

\bibitem[ABKP00]{MR2002g:57051}
J.~Amor{\'o}s, F.~Bogomolov, L.~Katzarkov, and T.~Pantev.
\newblock Symplectic {L}efschetz fibrations with arbitrary fundamental groups.
\newblock {\em J. Differential Geom.}, 54(3):489--545, 2000.
\newblock With an appendix by Ivan Smith.

\bibitem[EM02]{MR2003g:53164}
Y.~Eliashberg and N.~Mishachev.
\newblock {\em Introduction to the {$h$}-principle}, volume~48 of {\em Graduate
  Studies in Mathematics}.
\newblock American Mathematical Society, Providence, RI, 2002.

\bibitem[GLS96]{MR98d:58074}
Victor Guillemin, Eugene Lerman, and Shlomo Sternberg.
\newblock {\em Symplectic fibrations and multiplicity diagrams}.
\newblock Cambridge University Press, Cambridge, 1996.

\bibitem[GM03]{MR1950475}
Sergei~I. Gelfand and Yuri~I. Manin.
\newblock {\em Methods of homological algebra}.
\newblock Springer Monographs in Mathematics. Springer-Verlag, Berlin, second
  edition, 2003.

\bibitem[Gom95]{MR96j:57025}
Robert~E. Gompf.
\newblock A new construction of symplectic manifolds.
\newblock {\em Ann. of Math. (2)}, 142(3):527--595, 1995.

\bibitem[Gro86]{MR90a:58201}
Mikhael Gromov.
\newblock {\em Partial differential relations}, volume~9 of {\em Ergebnisse der
  Mathematik und ihrer Grenzgebiete (3) [Results in Mathematics and Related
  Areas (3)]}.
\newblock Springer-Verlag, Berlin, 1986.

\bibitem[HV04]{haller}
S.~Haller and C~Vizman.
\newblock Non-linear grassmannians as coadjoint orbits.
\newblock {\em Math. Ann.}, 329(4):771--785, 2004.

\bibitem[JK]{JK}
Tadeusz Januszkiewicz and Jarek K{\c e}dra.
\newblock Characteristic classes of smooth fibrations.
\newblock {\em math.SG/0209288}.

\bibitem[KM05]{MR2115670}
Jarek K{\c{e}}dra and Dusa McDuff.
\newblock Homotopy properties of {H}amiltonian group actions.
\newblock {\em Geom. Topol.}, 9:121--162 (electronic), 2005.

\bibitem[KN96a]{MR97c:53001a}
Shoshichi Kobayashi and Katsumi Nomizu.
\newblock {\em Foundations of differential geometry. {V}ol. {I}}.
\newblock Wiley Classics Library. John Wiley \& Sons Inc., New York, 1996.

\bibitem[KN96b]{MR97c:53001b}
Shoshichi Kobayashi and Katsumi Nomizu.
\newblock {\em Foundations of differential geometry. {V}ol. {II}}.
\newblock Wiley Classics Library. John Wiley \& Sons Inc., New York, 1996.

\bibitem[McD06]{dusa}
Dusa McDuff.
\newblock Enlarging the hamiltonian group.
\newblock {\em Journal of Symplectic Geometry}, 3(4):1--50, 2006.

\bibitem[MM65]{MR0174052}
John~W. Milnor and John~C. Moore.
\newblock On the structure of {H}opf algebras.
\newblock {\em Ann. of Math. (2)}, 81:211--264, 1965.

\bibitem[MS98]{MR2000g:53098}
Dusa McDuff and Dietmar Salamon.
\newblock {\em Introduction to symplectic topology}.
\newblock Oxford Mathematical Monographs. The Clarendon Press Oxford University
  Press, New York, second edition, 1998.

\bibitem[NR61]{MR0133772}
M.~S. Narasimhan and S.~Ramanan.
\newblock Existence of universal connections.
\newblock {\em Amer. J. Math.}, 83:563--572, 1961.

\bibitem[Seg70]{MR43:6292}
Graeme Segal.
\newblock Cohomology of topological groups.
\newblock In {\em Symposia Mathematica, Vol. IV (INDAM, Rome, 1968/69)}, pages
  377--387. Academic Press, London, 1970.

\bibitem[Spa66]{MR35:1007}
Edwin~H. Spanier.
\newblock {\em Algebraic topology}.
\newblock McGraw-Hill Book Co., New York, 1966.

\bibitem[tD87]{MR89c:57048}
Tammo tom Dieck.
\newblock {\em Transformation groups}, volume~8 of {\em de Gruyter Studies in
  Mathematics}.
\newblock Walter de Gruyter \& Co., Berlin, 1987.

\end{thebibliography}
\bibliographystyle{alpha}
\parindent=0pt

\end{document}